\documentclass{jhrs}
\newtheorem{theorem}{Theorem}
\newtheorem{definition}[theorem]{Definition}
\newtheorem{corollary}[theorem]{Corollary}

\newtheorem{example}[theorem]{Example}
\newtheorem{lemma}[theorem]{Lemma}

\newtheorem{proposition}[theorem]{Proposition}

\def\hnida#1{
{
\unitlength=.5pt
\begin{picture}(40.00,40.00)(0.00,0.00)
\put(38,20){\makebox(0.00,0.00)[l]{$#1$}}
\put(24.00,2.00){\makebox(0.00,0.00){$...$}}
\put(24.00,38.00){\makebox(0.00,0.00){$...$}}
\put(20.00,20.00){\makebox(0.00,0.00){$\bullet$}}
\put(40.00,0.00){\line(-1,1){20.00}}
\put(10.00,0.00){\line(1,2){10.00}}
\put(0.00,0.00){\line(1,1){20.00}}
\put(20.00,20.00){\line(1,1){20.00}}
\put(20.00,20.00){\line(-1,2){10.00}}
\put(20.00,20.00){\line(-1,1){20.00}}
\end{picture}}
}

\def\cal{\mathcal} \def\dd{{\mathbf s}} \def\gr{{\rm gr}}
\def\zn#1{{|\theta|(|\phi_1| + \cdots + |\phi_{#1}|)}}
\def\ai{$A_\infty$}\def\sfDB{{\sf DB}} \def\sfGS{{\sf GS}}
\def\calO{{\cal O}}\def\calS{{\cal S}}
\def\ggg#1#2{{\it Dgr\/}(#1,#2)}\def\edge{{\it edge}} \def\vert{{\it Vert}}
\def\sigmaspan{{\it Span}_{\mbox{\scriptsize $\Sigma$-$\Sigma$}}}
\def\Span{{\it Span\/}} \def\psfb{{\partial_{\sf B}}}
\def\whb{\widehat \beta} \def\Alg#1#2{{#1 \hskip -.1em  {\it Alg}(#2)}}
\def\lv{{\left(\rule{0em}{1em}\right. \hskip -.1em}}
\def\rv{{\hskip -.1em \left.\rule{0em}{1em}\right)}}
\def\ACA{{\it ACA}} \def\CA{{\it CA}}\def\AC{{\it AC}} 
\def\CAC{{\it CAC}} \def\calF{{\cal F}}
\def\sqldots{{\mbox{$. \hskip -.3pt . \hskip -.3pt .$}}}
\def\modU{{\rm U}} \def\copr{{\hskip 0.05em * \hskip .1em}}
\def\E{{\it Ex\/}} \def\w{\widetilde} \def\DER{\Der(\sfM,\sfM \copr \calE)}
 \def\B{{i_\sfM}_*} \def\C{\whb_*}
\def\sigmabimodmaps{\Lin_{\mbox{\scriptsize$\Sigma$-$\Sigma$}}}
\def\desusp{\downarrow\!} \def\In{{\it in\/}} \def\CB{C_\sfB}
\def\CBB{\it CB}\def\CcP{C_{\cal P}} \def\bpd{{\overline{\partial}_D}}

\def\susp{\uparrow\!} \def\CP{C_\sfP} \def\Out{{\it out\/}}
\def\squeezedldots{{\mbox{$. \hskip -1pt . \hskip -1pt .$}}}
\def\angles#1{{\langle #1 \rangle}} \def\sb#1{\langle #1 \rangle}
\def\Vert{{\it Vert\/}} \def\calP{{\cal P}} \def\biar{{\it biar\/}}
\def\GS{{\it GS}} \def\ot{\otimes} \def\otexp#1#2{#1^{\ot #2}}
\def\PROP{{\sc prop}} \def\calE{{\cal E}} \def\bfk{{\bf k}}
\def\Der{\mbox{\it Der\/}} \def\Ass{\mbox{${\cal A}${\it ss\/}}}
\def\ssAss{\mbox{\scriptsize${\cal A}${\it ss\/}}}

\def\Lin{{\mbox {\it Lin\/}}} \def\di{\diamond}
\def\id{{1 \!\! 1}} 
\def\iden{{\it id}}\def\li{$L_\infty$}
\def\deesusp#1#2{\mbox{$\downarrow^{#1} \hskip -.2em #2$}}
\def\rada#1#2{#1,\dots,#2}
\def\Rada#1#2#3{#1_{#2},\ldots,#1_{#3}}
\def\sqRada#1#2#3{#1_{#2},\mbox{$. \hskip -.2pt . \hskip -.2pt .$},#1_{#3}} 
\def\pa{\partial}
\def\otRada#1#2#3{#1_{#2}\otimes \cdots \otimes #1_{#3}}
\def\End{\hbox{${\sf End}$}} \def\freePROP{{\sf F}}
\def\Krtecek#1{\Upsilon_{s,{\bf v},{\bf t},#1}}
\def\Kr{\Krtecek}\def\brd{{\overline{\rho}_D}}
\def\colim#1{\mathop{{\rm colim}}%
             \limits_{\rule{0em}{.8em}\mbox{\scriptsize $#1$}}}
\def\UGR{\tt UGr}\def\UGr#1#2{{\UGR}(#1,#2)}
\def\sfB{{\sf B}}\def\sfM{{\sf M}} \def\sfP{{\sf P}} 
\def\CGSbi#1#2{{C^{#1,#2}_{GS}(B;B)}} \def\hPROP{$\frac12${\sc prop}}
 \def\papert{\partial_{\it pert\/}}
\def\SFM{\sfM}\def\sfDM{{\sf DM}}

\def\dvajednadec#1{{
\unitlength=.07em
\begin{picture}(24.00,20.00)(-2.00,0.00)
\put(10.00,10.00){\line(0,-1){10.00}}
\put(10.00,17.00){\makebox(0,0)[b]{\scriptsize $#1$}}
\put(10.00,10.00){\makebox(0,0){$\bullet$}}
\bezier{50}(10.00,10.00)(15.00,15.00)(20.00,20.00)
\bezier{50}(0.00,20.00)(5.00,15.00)(10.00,10.00)
\end{picture}}
}

\def\jednadvadec#1{{
\unitlength=.07em
\begin{picture}(24.00,20.00)(-2.00,0.00)
\bezier{50}(10.00,10.00)(15.00,5.00)(20.00,0.00)
\bezier{50}(10.00,10.00)(5.00,5.00)(0.00,0.00)
\put(12.00,15.00){\makebox(0,0)[lb]{\scriptsize $#1$}}
\put(10.00,5.00){\makebox(0,0)[b]{$\bullet$}}
\put(10.00,20.00){\line(0,-1){10.00}}
\end{picture}}
}

\def\bZbbZdec #1#2{
{
\unitlength=.05em
\begin{picture}(48.00,30.00)(-4,0.00)
\put(20.00,20.00){\line(1,-1){20}}
\put(20.00,20.00){\line(-1,-1){20}}
\bezier{50}(30.00,10.00)(25.00,5.00)(20.00,0.00)
\put(20.00,30.00){\line(0,-1){10.00}}
\put(23.00,26.00){\makebox(0,0)[lb]{\scriptsize $#1$}}
\put(31.00,18.00){\makebox(0,0)[lb]{\scriptsize $#2$}}
\put(20.00,20.00){\makebox(0,0){$\bullet$}}
\put(30.00,10.00){\makebox(0,0){$\bullet$}}
\end{picture}}
}

\def\ZbbZbdec#1#2{{
\unitlength=.05em
\begin{picture}(48.00,30.00)(-4,0.00)
\put(20.00,20.00){\line(1,-1){20}}
\put(20.00,20.00){\line(-1,-1){20}}
\bezier{50}(10.00,10.00)(15.00,5.00)(20.00,0.00)
\put(20.00,30.00){\line(0,-1){10.00}}
\put(23.00,26.00){\makebox(0,0)[lb]{\scriptsize $#1$}}
\put(10.00,18.00){\makebox(0,0)[rb]{\scriptsize $#2$}}
\put(20.00,20.00){\makebox(0,0){$\bullet$}}
\put(10.00,10.00){\makebox(0,0){$\bullet$}}
\end{picture}}}

\def\bVbbVdec #1#2{
{
\unitlength=.05em
\begin{picture}(48.00,30.00)(-4,-28.00)
\put(20.00,-20.00){\line(1,1){20}}
\put(20.00,-20.00){\line(-1,1){20}}
\bezier{50}(30.00,-10.00)(25.00,-5.00)(20.00,0.00)
\put(20.00,-30.00){\line(0,1){10.00}}
\put(23.00,-26.00){\makebox(0,0)[lt]{\scriptsize $#1$}}
\put(31.00,-18.00){\makebox(0,0)[lt]{\scriptsize $#2$}}
\put(20.00,-20.00){\makebox(0,0){$\bullet$}}
\put(30.00,-10.00){\makebox(0,0){$\bullet$}}
\end{picture}}
}

\def\VbbVbdec#1#2{{
\unitlength=.05em
\begin{picture}(48.00,30.00)(-4,-28.00)
\put(20.00,-20.00){\line(1,1){20}}
\put(20.00,-20.00){\line(-1,1){20}}
\bezier{50}(10.00,-10.00)(15.00,-5.00)(20.00,0.00)
\put(20.00,-30.00){\line(0,1){10.00}}
\put(23.00,-26.00){\makebox(0,0)[lt]{\scriptsize $#1$}}
\put(10.00,-18.00){\makebox(0,0)[rt]{\scriptsize $#2$}}
\put(20.00,-20.00){\makebox(0,0){$\bullet$}}
\put(10.00,-10.00){\makebox(0,0){$\bullet$}}
\end{picture}}}

\def\jednatridec#1{{
\unitlength=.07em
\begin{picture}(24.00,20.00)(-2.00,0.00)
\bezier{100}(10.00,10.00)(15.00,5.00)(20.00,0.00)
\bezier{100}(10.00,10.00)(5.00,5.00)(0.00,0.00)
\put(10.00,20.00){\line(0,-1){20.00}}
\put(13.00,15.00){\makebox(0,0)[lb]{\scriptsize $#1$}}
\put(10.00,10.00){\makebox(0,0){$\bullet$}}
\end{picture}}
}

\def\dvojiteypsilondec#1#2{{
\unitlength=.07em
\begin{picture}(24.00,30.00)(0.00,3.00)
\put(10.00,20.00){\line(0,-1){10.00}}
\bezier{70}(10.00,10.00)(15,5)(20.00,0.00)
\bezier{70}(10.00,10.00)(5,5)(0.00,0.00)
\bezier{70}(10.00,20.00)(15,25)(20.00,30.00)
\bezier{70}(0.00,30.00)(5,25)(10.00,20.00)
\put(10.00,20.00){\makebox(0,0){$\bullet$}}
\put(10.00,10.00){\makebox(0,0){$\bullet$}}
\put(10.00,26.00){\makebox(0,0)[b]{\scriptsize $#1$}}
\put(10.00,4.00){\makebox(0,0)[t]{\scriptsize $#2$}}
\end{picture}}}

\def\dvojiteypsilon{{
\unitlength=.3pt
\begin{picture}(24.00,30.00)(0.00,3.00)
\put(10.00,20.00){\line(0,-1){10.00}}
\bezier{20}(10.00,10.00)(15,5)(20.00,0.00)
\bezier{20}(10.00,10.00)(5,5)(0.00,0.00)
\bezier{20}(10.00,20.00)(15,25)(20.00,30.00)
\bezier{20}(0.00,30.00)(5,25)(10.00,20.00)
\end{picture}}}

\def\zeroatwo#1{%
{
\unitlength=.5pt
\begin{picture}(36.00,30.00)(-3.00,10.00)
\put(10.00,20.00){\makebox(0.00,0.00)[l]{$#1$}}
\put(20.00,10.00){\line(0,-1){5.00}}
\put(10.00,10.00){\line(0,-1){5.00}}
\put(30.00,30.00){\line(-1,0){30.00}}
\put(30.00,10.00){\line(0,1){20.00}}
\put(0.00,10.00){\line(1,0){30.00}}
\put(0.00,30.00){\line(0,-1){20.00}}
\end{picture}}
}

\def\twoazero#1{%
{
\unitlength=.5pt
\begin{picture}(36.00,30.00)(-3.00,0.00)
\put(20.00,25.00){\line(0,-1){5.00}}
\put(10.00,25.00){\line(0,-1){5.00}}
\put(10.00,10.00){\makebox(0.00,0.00)[l]{$#1$}}
\put(30.00,20.00){\line(-1,0){30.00}}
\put(30.00,0.00){\line(0,1){20.00}}
\put(0.00,0.00){\line(1,0){30.00}}
\put(0.00,20.00){\line(0,-1){20.00}}
\end{picture}}
}

\def\zeroathree#1{
{
\unitlength=.5pt
\begin{picture}(48.00,30.00)(-4.00,10.00)
\put(20.00,20.00){\makebox(0.00,0.00){$#1$}}
\put(0.00,30.00){\line(0,-1){20.00}}
\put(40.00,30.00){\line(-1,0){40.00}}
\put(40.00,10.00){\line(0,1){20.00}}
\put(0.00,10.00){\line(1,0){40.00}}
\put(30.00,10.00){\line(0,-1){5.00}}
\put(20.00,10.00){\line(0,-1){5.00}}
\put(10.00,10.00){\line(0,-1){5.00}}
\end{picture}}
}

\def\jednadva{{
\unitlength=.4pt
\begin{picture}(24.00,20.00)(-2.00,0.00)
\bezier{20}(10.00,10.00)(15.00,5.00)(20.00,0.00)
\bezier{20}(10.00,10.00)(5.00,5.00)(0.00,0.00)
\put(10.00,20.00){\line(0,-1){10.00}}
\end{picture}}
}

\def\jednactyri{{
\unitlength=.05pt
\begin{picture}(176.00,160.00)(-8.00,0.00)
\put(80.00,100.00){\line(0,1){60.00}}
\bezier{20}(80.00,80.00)(100.00,30.00)(110.00,0.00)
\bezier{20}(80.00,80.00)(60.00,30.00)(50.00,0.00)
\bezier{20}(80.00,80.00)(120.00,40.00)(160.00,0.00)
\bezier{20}(80.00,80.00)(40.00,40.00)(0.00,0.00)
\put(80.00,100.00){\line(0,-1){20.00}}
\end{picture}}
}

\def\ctyrijedna{{
\unitlength=.05pt
\begin{picture}(176.00,160.00)(-8.00,-160.00)
\put(80.00,-100.00){\line(0,-1){60.00}}
\bezier{20}(80.00,-80.00)(100.00,-30.00)(110.00,0.00)
\bezier{20}(80.00,-80.00)(60.00,-30.00)(50.00,0.00)
\bezier{20}(80.00,-80.00)(120.00,-40.00)(160.00,0.00)
\bezier{20}(80.00,-80.00)(40.00,-40.00)(0.00,0.00)
\put(80.00,-80.00){\line(0,-1){40.00}}
\end{picture}}
}

\def\dvajedna{{
\unitlength=.4pt
\begin{picture}(24.00,20.00)(-2.00,0.00)
\put(10.00,10.00){\line(0,-1){10.00}}
\bezier{20}(10.00,10.00)(15.00,15.00)(20.00,20.00)
\bezier{20}(0.00,20.00)(5.00,15.00)(10.00,10.00)
\end{picture}}
}

\def\dvadva{{
\unitlength=.8pt
\begin{picture}(12.00,10.00)(-1.00,0.00)
\bezier{30}(0.00,0.00)(5.00,5.00)(10.00,10.00)
\bezier{30}(0.00,10.00)(5.00,5.00)(10.00,0.00)
\end{picture}}
}

\def\jednatri{{
\unitlength=.4pt
\begin{picture}(24.00,20.00)(-2.00,0.00)
\bezier{20}(10.00,10.00)(15.00,5.00)(20.00,0.00)
\bezier{20}(10.00,10.00)(5.00,5.00)(0.00,0.00)
\put(10.00,20.00){\line(0,-1){20.00}}
\end{picture}}
}

\def\trijedna{{
\unitlength=.4pt
\begin{picture}(24.00,20.00)(-2.00,-20.00)
\bezier{20}(10.00,-10.00)(15.00,-5.00)(20.00,0.00)
\bezier{20}(10.00,-10.00)(5.00,-5.00)(0.00,0.00)
\put(10.00,-20.00){\line(0,1){20.00}}
\end{picture}}
}

\def\dvatri{{
\unitlength=0.4pt
\begin{picture}(24.00,20.00)(-2.00,0.00)
\put(10.00,10.00){\line(0,-1){10.00}}
\bezier{30}(0.00,0.00)(10.00,10.00)(20.00,20.00)
\bezier{30}(0.00,20.00)(10.00,10.00)(20.00,0.00)
\end{picture}}
}

\def\tridva{{
\unitlength=.4pt
\begin{picture}(24.00,20.00)(-2.00,-20.00)
\put(10.00,-10.00){\line(0,1){10.00}}
\bezier{30}(0.00,0.00)(10.00,-10.00)(20.00,-20.00)
\bezier{30}(0.00,-20.00)(10.00,-10.00)(20.00,0.00)
\end{picture}}
}

\def\dvacarkatri{{
\unitlength=.2pt
\begin{picture}(40.00,50.00)(0.00,0.00)
\put(20.00,30.00){\line(0,-1){10.00}}
\put(20.00,0.00){\line(0,1){20}}
\bezier{20}(20.00,20.00)(30.00,10.00)(40.00,0.00)
\bezier{20}(20.00,20.00)(10.00,10.00)(0.00,0.00)
\bezier{20}(20.00,30.00)(30.00,40.00)(40.00,50.00)
\bezier{20}(0.00,50.00)(10.00,40.00)(20.00,30.00)
\end{picture}}
}

\def\gen#1#2{
\if #11
    \if #22 \jednadva \else \fi
\else
\fi
\if #12
    \if #22 \dvadva \else \fi
\else
\fi
\if #12
    \if #21 \dvajedna \else \fi
\else
\fi
\if #13
    \if #22 \tridva \else \fi
\else
\fi
\if #13
    \if #21 \trijedna \else \fi
\else
\fi
\if #12
    \if #23 \dvatri \else \fi
\else
\fi
\if #11
    \if #23 \jednatri \else \fi
\else
\fi
\if #11
    \if #24 \jednactyri \else \fi
\fi
\if #14
    \if #21 \ctyrijedna \else \fi
\fi
}

\def\bZbbZ{
{
\unitlength=.27pt
\begin{picture}(48.00,30.00)(-4,0.00)
\bezier{34}(20.00,20.00)(30.00,10.00)(40.00,0.00)
\bezier{34}(20.00,20.00)(10.00,10.00)(0.00,0.00)
\bezier{20}(30.00,10.00)(25.00,5.00)(20.00,0.00)
\put(20.00,30.00){\line(0,-1){10.00}}
\end{picture}}
}

\def\ZbbZb{{
\unitlength=.27pt
\begin{picture}(48.00,30.00)(-4,0.00)
\bezier{34}(20.00,20.00)(30.00,10.00)(40.00,0.00)
\bezier{34}(20.00,20.00)(10.00,10.00)(0.00,0.00)
\bezier{20}(10.00,10.00)(15.00,5.00)(20.00,0.00)
\put(20.00,30.00){\line(0,-1){10.00}}
\end{picture}}}

\def\ZbbZbb{{
\unitlength=.07pt
\begin{picture}(192.00,120.00)(-16.00,0.00)
\bezier{20}(20.00,20.00)(30.00,10.00)(40.00,0.00)
\bezier{30}(80.00,80.00)(120.00,40.00)(160.00,0.00)
\bezier{30}(80.00,80.00)(40.00,40.00)(0.00,0.00)
\put(80.00,120.00){\line(0,-1){120.00}}
\end{picture}}
}

\def\bbZbbZ{{
\unitlength=0.07pt
\begin{picture}(192.00,120.00)(-16.00,0.00)
\bezier{20}(140.00,20.00)(130.00,10.00)(120.00,0.00)
\bezier{30}(80.00,80.00)(120.00,40.00)(160.00,0.00)
\bezier{30}(80.00,80.00)(40.00,40.00)(0.00,0.00)
\put(80.00,120.00){\line(0,-1){120.00}}
\end{picture}}
}

\def\bZbbZb{{
\unitlength=.07pt
\begin{picture}(192.00,120.00)(-16.00,0.00)
\bezier{20}(80.00,20.00)(90.00,10.00)(100.00,0.00)
\bezier{20}(80.00,20.00)(70.00,10.00)(60.00,0.00)
\put(80.00,40.00){\line(0,-1){20.00}}
\put(80.00,40.00){\line(0,1){0.00}}
\put(80.00,80.00){\line(0,-1){40.00}}
\put(80.00,120.00){\line(0,-1){40.00}}
\bezier{30}(80.00,80.00)(120.00,40.00)(160.00,0.00)
\bezier{30}(80.00,80.00)(40.00,40.00)(0.00,0.00)
\end{picture}}
}

\def\ZbbbZb{{
\unitlength=.07pt
\begin{picture}(192.00,120.00)(-16.00,0.00)
\put(40.00,0.00){\line(0,1){20.00}}
\put(40.00,40.00){\line(0,-1){40.00}}
\bezier{20}(40.00,40.00)(60.00,20.00)(80.00,0.00)
\put(80.00,40.00){\line(0,1){0.00}}
\put(80.00,120.00){\line(0,-1){40.00}}
\bezier{30}(80.00,80.00)(120.00,40.00)(160.00,0.00)
\bezier{30}(80.00,80.00)(40.00,40.00)(0.00,0.00)
\end{picture}}
}

\def\bZbbbZ{{
\unitlength=0.07pt
\begin{picture}(192.00,120.00)(-16.00,0.00)
\put(120.00,40.00){\line(0,-1){40.00}}
\bezier{20}(120.00,40.00)(100.00,20.00)(80.00,0.00)
\put(80.00,40.00){\line(0,1){0.00}}
\put(80.00,120.00){\line(0,-1){40.00}}
\bezier{30}(80.00,80.00)(120.00,40.00)(160.00,0.00)
\bezier{30}(80.00,80.00)(40.00,40.00)(0.00,0.00)
\end{picture}}
}

\def\ZvvZv{{
\unitlength=.27pt
\begin{picture}(48.00,30.00)(-4,-30.00)
\bezier{30}(20.00,-20.00)(30.00,-10.00)(40.00,0.00)
\bezier{30}(20.00,-20.00)(10.00,-10.00)(0.00,0.00)
\bezier{20}(10.00,-10.00)(15.00,-5.00)(20.00,0.00)
\put(20.00,-30.00){\line(0,1){10.00}}
\end{picture}}}

\def\vZvvZ{{
\unitlength=.27pt
\begin{picture}(48.00,30.00)(-4,-30.00)
\bezier{34}(20.00,-20.00)(30.00,-10.00)(40.00,0.00)
\bezier{34}(20.00,-20.00)(10.00,-10.00)(0.00,0.00)
\bezier{20}(30.00,-10.00)(25.00,-5.00)(20.00,0.00)
\put(20.00,-30.00){\line(0,1){10.00}}
\end{picture}}
}

\def\vZvvZdva{{
\unitlength=.2pt
\begin{picture}(48.00,40.00)(-4.00,0.00)
\bezier{10}(30.00,30.00)(25.00,35.00)(20.00,40.00)
\bezier{34}(0.00,0.00)(20.00,20.00)(40.00,40.00)
\bezier{34}(0.00,40.00)(20.00,20.00)(40.00,0.00)
\end{picture}}}

\def\dvabZbbZ{{
\unitlength=.2pt
\begin{picture}(48.00,40.00)(-4.00,-40.00)
\bezier{10}(30.00,-30.00)(25.00,-35.00)(20.00,-40.00)
\bezier{34}(0.00,0.00)(20.00,-20.00)(40.00,-40.00)
\bezier{34}(0.00,-40.00)(20.00,-20.00)(40.00,0.00)
\end{picture}}}

\def\ZvvZvdva{{
\unitlength=.2pt
\begin{picture}(48.00,40.00)(-4.00,0.00)
\bezier{10}(10.00,30.00)(15.00,35.00)(20.00,40.00)
\bezier{34}(0.00,0.00)(20.00,20.00)(40.00,40.00)
\bezier{34}(0.00,40.00)(20.00,20.00)(40.00,0.00)
\end{picture}}}

\def\dvaZbbZb{{
\unitlength=.2pt
\begin{picture}(48.00,40.00)(-4.00,-40.00)
\bezier{10}(10.00,-30.00)(15.00,-35.00)(20.00,-40.00)
\bezier{34}(0.00,0.00)(20.00,-20.00)(40.00,-40.00)
\bezier{34}(0.00,-40.00)(20.00,-20.00)(40.00,0.00)
\end{picture}}}

\def\dvojicezob#1#2#3#4{{
{
\unitlength=.58pt
\thinlines
\begin{picture}(110.00,20.0)(0.00,14.00)
\put(49.50,0.00){\makebox(0.00,0.00)[t]{\scriptsize $#4$}}
\put(49.50,41.50){\makebox(0.00,0.00)[b]{\scriptsize $#3$}}
\put(97.00,20.00){\makebox(0.00,0.00)[lc]{$#2$}}
\put(5.00,20.00){\makebox(0.00,0.00)[rc]{$#1$}}
\put(23.00,7.5){\vector(-2,1){10.00}}
\put(79.00,31.00){\vector(2,-1){10.00}}
\bezier{100}(16.0,10.50)(49.50,0.50)(83.50,12.50)
\bezier{100}(16.0,30.50)(49.50,40.50)(83.50,28.50)
\end{picture}}
}}

\begin{document}
\title{Intrinsic brackets and the $L_\infty$-deformation theory of bialgebras}
\author{Martin Markl}

\email{markl@math.cas.cz}

\address{Math. Inst. of the Academy, {\v Z}itn{\'a} 25,
         115 67 Prague 1, The Czech Republic}

\thanks{The author was supported by the grant GA \v CR 201/02/1390 and by
   the Academy of Sciences of the Czech Republic,
   Institutional Research Plan No.~AV0Z10190503.}


\bibliographystyle{plain}
\hyphenation{coho-mo-logy}

\begin{abstract}
We show that there exists a Lie bracket on the cohomology of any
type of (bi)algebras over an operad or a \PROP, induced by an
$L_\infty$-structure on the defining cochain complex, such that the
associated $L_\infty$-master equation captures deformations.

This in particular implies the existence of
a Lie bracket on the Gerstenhaber-Schack
cohomology~\cite{gerstenhaber-schack:Proc.Nat.Acad.Sci.USA90} of a
bialgebra that
extends the classical intrinsic bracket~\cite{gerstenhaber:AM63} on
the Hochschild cohomology, giving an affirmative answer to an old
question about the existence of such a~bracket. We also explain how
the results of~\cite{umble-saneblidze:KK} provide explicit formulas
for this bracket.
\end{abstract}

\maketitle

\subsection*{Conventions}  
We assume a certain familiarity with operads and \PROP{s},
see~\cite{markl:ba,markl:handbook,markl-shnider-stasheff:book,mv,vallette:thesis}. 
The reader who wishes only to know how the intrinsic bracket on the
Gerstenhaber-Schack cohomology looks might proceed directly to
Section~\ref{GS} which is almost independent on the rest of the paper
and contains explicit calculations. We also assume some
knowledge of the concept of strongly homotopy Lie algebras (also
called $L_\infty$-algebras),
see~\cite{lada-markl:CommAlg95,lada-stasheff:IJTP93}.

We will make no distinction between an
operad $\calP$ and the
\PROP\ $\sfP$ generated by this operad. This means that
for us operads are particular cases of \PROP{s}.
As usual, bialgebra will mean a Hopf algebra without
(co)unit and antipode. To distinguish these bialgebras from other
types of ``bialgebras'' we will sometimes call them also
$\Ass$-bialgebras.

All algebraic objects will be defined over a fixed field $\bfk$ of
characteristic zero although, surprisingly, our constructions related
to $\Ass$-bialgebras make sense over the integers.

\section{Introduction and main results}
\label{intro}

We show that the cohomology of (bi)algebras always carries a Lie
bracket (which we call the {\em intrinsic bracket\/}), induced by an
$L_\infty$-structure on the corresponding cochain complex.  We also
discuss the master equation related to this $L_\infty$-structure.

By a {\em (bi)algebra\/} we mean an algebra over a certain $\bfk$-linear
\PROP\ $\sfP$.
Therefore a (bi)algebra is given by a
homomorphism of \PROP{s} $\alpha : \sfP \to \End_V$, where $\End_V$
denotes the endomorphism \PROP\ of a $\bfk$-vector space $V$. Observe that
this notion encompasses not only ``classical'' algebras (associative,
commutative associative, Lie, \&c.) but also various types of bialgebras
($\Ass$-bialgebras, Lie bialgebras, infinitesimal bialgebras, \&c.).

Let us recall that a {\em minimal model\/} of a $\bfk$-linear \PROP\
$\sfP$ is a differential (non-negatively) graded {\bfk}-linear
{\sc prop} $({\sf M},\pa)$ together with a homology isomorphism
\[
(\sfP,0) \stackrel{\rho}{\longleftarrow} ({\sf M},\pa)
\]
such that (i) the \PROP\ ${\sf M}$ is free and
(ii) the image of the degree $-1$ differential
$\pa$ consists of decomposable elements of ${\sf M}$ (the
minimality condition), see~\cite{markl:ba} for details. It is not our
aim to discuss in this paper the existence and uniqueness of minimal
models, nor the methods how to construct such models
explicitly. Let us say only that for a large class of operads and
\PROP{s} these minimal models can be constructed using the 
Koszul duality~\cite{gan,ginzburg-kapranov:DMJ94,vallette:thesis}.

Let us emphasize that instead of a minimal model of $\sfP$ we may use
in the following constructions any cofibrant (in a suitable sense)
resolution of $\sfP$. But since explicit minimal models of $\sfP$
exist in all cases of interest we will stick to minimal models in this
note. This will simplify some technicalities.

Assume we are given a homomorphism $\alpha : \sfP \to \End_V$
describing a $\sfP$-algebra $B$. To define its cohomology, we need to
choose first a minimal model $\rho : (\sfM,\pa) \to (\sfP,0)$ of
$\sfP$.  The composition $\beta:= \alpha \circ \rho : \sfM \to \End_V$
makes $\End_V$ an $\sfM$-module (in the sense
of~\cite[page~203]{markl:JPAA96}), one may therefore consider the
graded vector space of derivations $\Der(\sfM,\End_V)$. For $\theta
\in \Der(\sfM,\End_V)$ define $\delta \theta : = \theta \circ \pa$. It
follows from the obvious fact that $\beta \circ \pa = 0$, implied by
the triviality of the differential in $\sfP$, that $\delta \theta$ is
again a derivation, so $\delta$ is a well-defined endomorphism of
$\Der(\sfM,\End_V)$ which clearly satisfies $\delta^2=0$. We conclude
that $\Der(\sfM,\End_V)$ is a non-positively graded vector space
equipped with a differential of degree~$-1$. Finally, let
\begin{equation}
\label{Joshua_v_Praze}
\CP^*(V;V) := \susp \Der(\sfM,\End_V)^{-*}
\end{equation}
be the suspension of the graded vector space $\Der(\sfM,\End_V)$ with
reversed degrees.  The differential $\delta$ induces on $\CP^*(V;V)$ a
degree $+1$ differential denoted by $\delta_\sfP$.  The {\em cohomology
of $B$ with coefficients in itself\/} is then defined by
\begin{equation}
\label{Tel-Aviv}
H^*_{\sfP}(B;B) := H(\CP^*(V;V),\delta_\sfP),
\end{equation}
see~\cite{markl:JPAA96,markl:ba}.

For ``classical'' algebras, the cochain complex $(\CP^*(V;V),\delta_\sfP)$
agrees with the ``standard'' constructions.  Thus, for associative
algebras,~(\ref{Tel-Aviv}) gives the Hochschild cohomology, for
associative commutative algebras the Harrison cohomology, for Lie
algebras the Chevalley-Eilenberg cohomology,~\&c.  More generally, for
algebras over a quadratic Koszul operad the above cohomology coincides
with the triple cohomology. Therefore nothing dramatically new happens
here.

This situation changes if we consider (bi)algebras over a general
\PROP\ $\sfP$. To our best knowledge, (\ref{Tel-Aviv}) is the only
definition of a cohomology of (bi)algebras over \PROP{s}. As we
argued in~\cite{markl:JPAA96}, it governs deformations of these
(bi)algebras.

Let $(\sfB,0) \leftarrow (\sfM_\sfB,\pa)$ be the minimal model of the
\PROP\ $\sfB$ for $\Ass$-bialgebras constructed
in~\cite{markl:ba,umble-saneblidze:KK} and $B$ a bialgebra
given by a homomorphism $\alpha : \sfB \to \End_V$. Then
$(\CB^*(V;V),\delta_\sfB)$ is isomorphic to the Gerstenhaber-Schack
cochain complex $(C^*_\GS(B;B),d_\GS)$ and~(\ref{Tel-Aviv}) coincides
with the Gerstenhaber-Schack
cohomology~\cite{gerstenhaber-schack:Proc.Nat.Acad.Sci.USA90},
\[
H^*_\GS(B;B) \cong H^*(\CB^*(V;V),\delta_\sfB),
\]
see Section~\ref{GS} for details.
The result announced in the Abstract follows from the following:

\begin{theorem}
\label{main}
Let $B$ be a (bi)algebra over a \PROP\ $\sfP$. Then there exist a
graded Lie algebra bracket on the cohomology $H^*_\sfP(B;B)$ induced
by a natural $L_\infty$-structure $(\delta_\sfP, l_2,l_3,\ldots)$ on
the defining complex $(\CP^*(V;V),\delta_\sfP)$.
\end{theorem}

We will see in Section~\ref{1bis} and also in Section~\ref{3} that the
 $L_\infty$-structure of Theorem~\ref{main} can be given by explicit
 formulas that involve the differential $\pa$ of the minimal model
 $\sfM$.  It will also be clear that this $L_\infty$-structure uses
 all the information about the minimal model $\sfM$ of $\sfP$ and that,
 vice versa, the minimal model $\sfM$ can be reconstructed from the
 knowledge of this $L_\infty$-structure. Therefore the brackets
 $l_2,l_3,\ldots$ can be understood as Massey products that detect the
 homotopy type of the \PROP~$\sfP$.

Since the minimal model $\sfM$ is, by definition, free on a
$\Sigma$-bimodule $E$, $\sfM = \freePROP(E)$, it is graded by the number
of generators. This means that $\freePROP(E) = \bigoplus_{k \geq 0}
\freePROP^k(E)$, where $\freePROP^k(E)$ is spanned by ``monomials'' composed of
exactly $k$ elements of $E$.  The minimality of $\pa$ is equivalent to
$\pa(E) \subset \freePROP^{\geq 2}(E)$.  The differential $\pa$ is called
{\em quadratic\/} if $\pa(E) \subset \freePROP^2(E)$.

\begin{proposition}
\label{main1}
If the differential of the minimal model $\sfM$ of $\sfP$ used in the
definition of the cohomology~(\ref{Tel-Aviv}) is quadratic, then the
higher brackets $l_3,l_4,\ldots$ of the $L_\infty$-structure vanish,
therefore $(\CP^*(V;V),\delta_\sfP)$ forms an ordinary dg Lie algebra with
the bracket $[-,-]:= l_2(-,-)$.
\end{proposition}

The minimal model of a quadratic Koszul operad $\calP$ is given by the
cobar construction on its quadratic dual $\calP^!$ and is therefore
quadratic.  Thus, for algebras over such an operad, the complex
$(C_\calP^*(V;V),\delta_\calP)$ is a Lie algebra whose bracket
coincides with the classical intrinsic bracket given by identifying
this complex with the space of coderivations of a certain cofree
nilpotent $\calP^!$-coalgebra,
see~\cite[Section~II.3.8]{markl-shnider-stasheff:book}.  The similar
observation is true also for various types of ``bialgebras'' defined
over quadratic (in a suitable sense) \PROP{s}, such as Lie
bialgebras~\cite{gan}, infinitesimal bialgebras~\cite{aguiar:04} and
$\frac12$bialgebras~\cite{markl:ba}.  In contrast, as we will see in
Section~\ref{GS}, the Gerstenhaber-Schack complex of an
$\Ass$-bialgebra carries a fully fledged $L_\infty$-algebra structure.

\vskip .3em
\noindent
{\bf Relation to previous results.}  As indicated in the above
paragraph, it is well-known that, for an algebra $B$ over a quadratic
Koszul operad $\calP$, the cochain complex
$(C_\calP^*(V;V),\delta_\calP)$ is a dg-Lie algebra with the structure
given by a generalization of Schlessinger-Stasheff's intrinsic
bracket~\cite{schlessinger-stasheff:JPAA85}. For algebras over a
general operad, an $L_\infty$-generalization of this structure was
obtained by van~der~Laan~\cite{laan-defo} as follows.

Van~der~Laan noticed that, for each homotopy cooperad (in an
appropriate sense) $E$ and for each operad $\calS$, the
$\Sigma$-module $\calS^E = \{\calS^E(n)\}_{n \geq 1}$, where
$\calS^E(n) := \Lin(E(n),\calS(n))$, is a homotopy operad (again in an
appropriate sense), which generalizes the convolution operad
of~\cite{berger-moerdijk:02}. He also proved that, for each homotopy
operad $\calO = \{\calO(n)\}_{n \geq 1}$, the ``total
space'' $\calO^* := \bigoplus_{* \geq 1}\calO(*+1)$ has an induced
$L_\infty$-structure which descents to an $L_\infty$-structure on the
symmetrization $\calO^*_\Sigma := \bigoplus_{* \geq
1}\calO(*+1)_{\Sigma_{*+1}}$.  Therefore, for $\calS$ and $E$ as above, the
graded vector space ${\calS^E}^*_\Sigma$ is a natural $L_\infty$-algebra.

On the other hand, let $(\freePROP(E),\pa) \to (\calP,0)$ be a minimal
model of $\calP$. Van~der~Laan observed that the differential $\pa$
makes the $\Sigma$-module of generators $E$ a homotopy cooperad and that,
for $\calS = \End_V$,
\begin{equation}
\label{zavody_ve_Dvore_zacaly_a_porad_prsi}
{\calS^E}^*_\Sigma \cong C_\calP^*(V;V).
\end{equation}
Combining the above facts, 
he concluded that $\CP^*(V;V)$ is a natural $L_\infty$-algebra
and proved that the map $\alpha : \calP \to \End_V$ defining the
$\calP$-algebra $B$ determines a Maurer-Cartan element $\kappa \in
C_\calP^1(V;V)$.  He then constructed the $L_\infty$-structure on
$(C_\calP^*(V;V),\delta_\calP)$ as the $\kappa$-twisting, in the sense
recalled in Section~\ref{sh}, of the $L_\infty$-algebra given
by the identification~(\ref{zavody_ve_Dvore_zacaly_a_porad_prsi}).

We were recently informed about an on-going 
work~\cite{merkulov-vallette} whose central statement
proves the existence of an $L_\infty$-structure on
the space of $\mathbb Z$-graded extended morphisms from a free dg \PROP\ to an
arbitrary \PROP\ from which Laan's arguments and their generalization
to \PROP{s} follow.

The methods of this article are independent of the above mentioned
papers.  While our approach is not very conceptual, it is
straightforward and immediately produces, from a given differential
in the minimal model, explicit formulas for the induced
$L_\infty$-structure.

\vskip .3em
\noindent
{\bf Relation to derived spaces of algebra structures.}  
As argued in~\cite[page~797]{fontanine-kapranov:JAMS}, for an operad
$\calP$ and a {\em finite-dimensional\/} vector space $W$, there
exists a scheme $\Alg {\calP}W$ parameterizing $\calP$-algebra
structures on $W$.  It is characterized by the property that for each
commutative dg-algebra $A$, morphisms ${\it Spec\/}(A) \to \Alg
{\calP}W$ are in bijection with $(A \ot_\bfk \calP)$-algebra
structures on the dg-$A$-module $A \ot_\bfk W$.  Let $\calF \to \calP$
be a free resolution of the operad $\calP$. Then $\Alg {\calF}W$ was
interpreted, in~\cite[Section~3.2]{fontanine-kapranov:JAMS}, as a
smooth dg-scheme in the sense of~\cite{fontanine-kapranov},
representing a right-derived space $R\hskip .1em \Alg {\calP}W$ of
$\calP$-actions on $W$ in a suitable derived category of dg-schemes.

It can be easily seen, using methods
of~\cite[Section~3.5]{fontanine-kapranov:JAMS}, that if $\rho: \calF
\to \calP$ is the minimal model of the operad $\calP$ and $\alpha :
\calP \to \End_W$ describes a $\calP$-algebra $B$ with the underlying
vector space $W$, then the components of the dg-tangent space at
$[\beta] \in \Alg {\calF}W$, $\beta := \rho \circ \alpha$, can be
described as
\[
T^n_{[\beta]} \Alg {\calF}W \cong \CcP^{n+1}(V;V),\ \mbox {for } n \geq 0
\]   
(we used a different degree convention than~\cite{fontanine-kapranov:JAMS}).
The existence of an $L_\infty$-structure on $\CcP^{*}(V;V)$ would then
follow from general properties of dg-schemes and is in fact
equivalent to specifying a local coordinate system at the smooth point
$[\beta]$ of $\Alg {\calF}W$~\cite[Proposition~2.5.8]{fontanine-kapranov}.

On the other hand, let $\sfM$ be a free dg-\PROP\ 
and $\beta : \sfM \to \calE$ a
\PROP\ homomorphism. Denote by $\calE_\beta$ the \PROP\ $\calE$
considered as an $\sfM$-module with the action induced by the
homomorphism~$\beta$. We will prove in Theorem~\ref{E} of
Section~\ref{2bis} that the desuspended space $\desusp
\Der(\SFM,\calE_\beta)$ of derivations has a natural
$L_\infty$-structure. By the definition~(\ref{Joshua_v_Praze}) of
$\CP^*(V,V)$, Theorem~\ref{main} follows from Theorem~\ref{E}
by taking $\sfM$ a minimal model of the
$\sfP$, $\calE: = \End_V$ and $\beta: \sfM \to \calE$ 
the composition $\alpha \circ
\rho$, where $\alpha : \sfP \to \End_V$ describes the algebra $B$ and
$\rho : \sfM \to \sfP$ is the map of the minimal~model.

In the light of~\cite[Proposition~2.5.8]{fontanine-kapranov},
Theorem~\ref{E} translates to the statement that for each free
dg-\PROP\ $\sfM$ and each homomorphism $\beta : \sfM \to \calE$, the space
$\desusp \Der(\sfM,\calE_\beta)$ forms a smooth dg-scheme.  The
derived scheme $R\hskip .1em \Alg {\calP}W$
of~\cite[page~797]{fontanine-kapranov:JAMS} is then, for $\sfP$ the
operad $\calP$ and $W$ a finite-dimensional vector space, the
specialization of this construction at the point represented by $\calE
= \End_W$ and $\beta = \alpha \circ \rho$.  It this sense, the results
of the present paper are meta-versions of constructions
in~\cite[Section~3.2]{fontanine-kapranov:JAMS} that completely avoid
all assumptions required by the
`classical' geometry, namely the fact 
that the target of the map $\beta$ is the endomorphism \PROP\ of a
finite-dimensional vector space. The present paper thus 
finishes the program to find a ``universal variety of
structure constants'' formulated in~\cite[page~197]{markl:JPAA96}. 

\vskip .3em

\noindent
{\bf The master equation.}  Let $A$ be a ``classical'' algebra over a
quadratic Koszul operad $\calP$ (associative, commutative associative, Lie,
\&c.), so that the cochain complex $(C_\calP^*(A;A),\delta_\calP)$ is a
graded dg-Lie algebra (Proposition~\ref{main1}). One usually shows that an
element $\kappa \in C_\calP^1(A,A)$, represented by a bilinear map (or
by a collection of bilinear maps), is a deformation of the 
$\calP$-algebra structure $A$
if and only if it solves the ``classical'' master equation $0 =
\delta_\calP(\kappa) +
\textstyle\frac 12 [\kappa,\kappa]$.

{}For a (bi)algebra $B$ over a general \PROP\ $\sfP$, the
cochain complex $(C^*_\sfP(B,B),d_\sfP)$ forms only an $L_\infty$-algebra,
but we will prove, in Section~\ref{sh},
that solutions $\kappa \in C_\sfP^1(B,B)$ of the
``quantum'' master equation
\begin{equation}
\label{master1}
0 = \delta_\sfP(\kappa) + \frac 1{2!}  l_2(\kappa,\kappa)
- \frac 1{3!}  l_3(\kappa,\kappa,\kappa) -
\frac 1{4!}  l_4(\kappa,\kappa,\kappa,\kappa) + \cdots
\end{equation}
are deformations of $B$. This means that the $L_\infty$-structure
of Theorem~\ref{main} represents an $L_\infty$-version of the Deligne
groupoid for deformations of $B$, see~\cite{hinich:IMRN97} for
the terminology.
We will see in Section~\ref{3} how this observation applies to
$\Ass$-bialgebras. Although the sum~(\ref{master1}) is infinite, we
will see that, in situations considered in this paper, it converges. 

\vskip .3em

\noindent
{\bf Acknowledgment:} The author would like to express his thanks to
the Mathematics Department of Bar Ilan University for a very
stimulating and pleasant atmosphere. It was during his visit of this
department, sponsored by the Israel Academy of Sciences, that the
present work started.  Also the suggestions and comments of Ezra
Getzler, Jim Stasheff, Bruno Vallette and Sasha Voronov regarding the
first draft of this paper were extremely
helpful. Proposition~\ref{podzim} and most of Section~\ref{sh} is
basically only a \PROP{ic} generalization of the material on
pages~21-23 of~\cite{laan-defo}.  My thanks are also due to
S.~Merkulov and B.~Vallette for sharing their on-going
work~\cite{merkulov-vallette} with me.
\vskip .3em

\noindent
{\bf Outline of the paper.}
In the following section we indicate the idea behind the
$L_\infty$-structure of Theorem~\ref{main}. A rigorous proof of this
theorem is then contained in Sections~\ref{2} and~\ref{2bis}. In short
Section~\ref{sh} we discuss master equations in $L_\infty$-algebras.
In the last section we show how constructions of this paper together with the
description~\cite[Eqn.~3.1]{umble-saneblidze:KK} of the minimal model of the
bialgebra \PROP\ give an explicit $L_\infty$-structure on the
Gerstenhaber-Schack complex $(C^*_\GS(B;B),d_\GS)$.


\section{The idea of the construction}
\label{1bis}

In this section we explain the idea behind the $L_\infty$-structure of
Theorem~\ref{main} and indicate why the $L_\infty$-axioms are
satisfied.  We believe that this section will help to
understand the concepts, but we do not aim to be rigorous here, see
also the remark at the end of this section. Formal constructions and
proofs based on the equivalence between symmetric brace algebras and
pre-Lie algebras~\cite{guin-oudom,lm:sb} are then given in
Sections~\ref{2} and~\ref{2bis}.

We need to review first some definitions and facts concerning \PROP{s} and
their derivations.  Given a \PROP\ $\sfP$ and a $\sfP$-module
$\modU$~\cite[p.~203]{markl:JPAA96}, then a {\em degree $d$ derivation\/}
$\theta : \sfP \to \modU$ is a map of $\Sigma$-bimodules $\theta
: \sfP \to \modU$ which is a degree $d$ derivation (in the
evident sense) with respect to both the horizontal and vertical
compositions in the \PROP\ $\sfP$ and the $\sfP$-module~$\modU$.

An equivalent definition is the following. For each $d$, the
$\sfP$-module structure on $\modU$ induces the obvious \PROP\
structure on the direct sum $\sfP \oplus \deesusp d{\modU}$ of
the $\Sigma$-bimodule $\sfP$ and the $d$-fold desuspension of the
$\Sigma$-bimodule $\modU$. A degree $d$ map $\theta : \sfP \to
\modU$ of $\Sigma$-bimodules is then a degree $d$ derivation if
and only if
\begin{equation}
\label{jitka3}
{\it id}_\sfP \oplus \deesusp d\theta :
\sfP \to \sfP \oplus \deesusp d{\modU}
\end{equation}
is a \PROP\ homomorphism.  The equivalence of the above two
definitions of derivations can be easily verified directly. 
We denote
by $\Der(\sfP,\modU)$ the graded vector space of derivations $\theta :
\sfP \to \modU$. 
If $\modU = \sfP$, we write simply $\Der(\sfP)$ instead of
$\Der(\sfP,\sfP)$.

\begin{proposition}
\label{JITa}
Let $\sfM = \freePROP(E)$ be the free \PROP\ generated by a
$\Sigma$-bimodule $E$ and $\modU$ an $\sfM$-module. Then there is
a canonical isomorphism
\begin{equation}
\label{jitka2} 
\Der(\sfM,\modU) \cong \sigmabimodmaps(E,\modU),
\end{equation}
given by restricting a derivation $\theta \in \Der(\sfM,\modU)$ onto
the space $E \subset \sfM$ of generators. In~(\ref{jitka2}),
$\sigmabimodmaps(-,-)$ denotes the space of linear
bi-equivariant maps of $\Sigma$-bimodules.
\end{proposition}

\noindent 
The {\it proof\/} follows from the
interpretation~(\ref{jitka3}) of derivations as homomorphisms and
the standard universal property of free \PROP{s}.%
\qed

Let us look more closely at the structure of the free \PROP\
$\freePROP(E)$ generated by a
$\Sigma$-bimodule~$E$. As~explained in~\cite{markl:handbook}, the
components of this \PROP\ are the colimit
\begin{equation}
\label{jsem_zvedav_jestli_se_ta_druha_Jitka_ozve}
\freePROP(E)(m,n) := \colim{{G \in \UGr mn}}{{E}}(G),\
m,n \geq 0,
\end{equation} 
taken over the category $\UGr mn$ of directed $(m,n)$-graphs without
directed cycles and their
isomorphisms. In~(\ref{jsem_zvedav_jestli_se_ta_druha_Jitka_ozve}),
$E(G)$ denotes the vector space of all decorations of vertices of $G$
by elements of $E$, see~\cite[Section~8]{markl:handbook} for precise
definitions. Therefore elements of the free \PROP\ $\freePROP(E)$ can
be represented by sums of $E$-decorated directed~graphs.

To simplify the exposition, we accept the convention that $\Gamma$
(with or without a subscript) will denote an $E$-decorated graph, and
$G$ (with or without a subscript) the underlying un-decorated graph.
If $\Gamma$ is such an $E$-decorated graph, we denote by
$e_v \in E$ the corresponding decoration of a vertex $v \in \Vert(G)$ 
of the underlying un-decorated graph. 

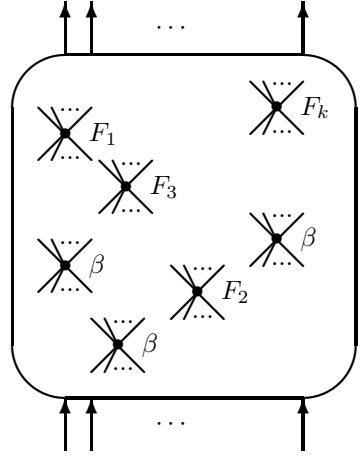
\begin{figure}[t]
\begin{center}
{
\unitlength=1pt
\begin{picture}(170.00,170.00)(0.00,0.00)
\thicklines
\put(90.00,10.00){\makebox(0.00,0.00){$\cdots$}}
\put(90.00,160.00){\makebox(0.00,0.00){$\cdots$}}
\put(130.00,80.00){\makebox(0.00,0.00){$\hnida {\beta}$}}
\put(70.00,40.00){\makebox(0.00,0.00){$\hnida {\beta}$}}
\put(50.00,70.00){\makebox(0.00,0.00){$\hnida {\beta}$}}
\put(130.00,130.00){\makebox(0.00,0.00){$\hnida {F_k}$}}
\put(100.00,60.00){\makebox(0.00,0.00){$\hnida {F_2}$}}
\put(73.00,100.00){\makebox(0.00,0.00){$\hnida {F_3}$}}
\put(50.00,120.00){\makebox(0.00,0.00){$\hnida {F_1}$}}
\put(140.00,150.00){\vector(0,1){20.00}}
\put(60.00,150.00){\vector(0,1){20.00}}
\put(50.00,150.00){\vector(0,1){20.00}}
\put(140.00,0.00){\vector(0,1){20.00}}
\put(60.00,0.00){\vector(0,1){20.00}}
\put(50.00,0.00){\vector(0,1){20.00}}
\put(95.00,85.00){\oval(130.00,130.00)}
\end{picture}}
\end{center}
\caption{\label{figu1} 
The $\End_V$-decorated graph $\Gamma_{\{\beta\}}^{\{\Rada
v1k\}}[\Rada F1k]$. Vertices labelled $F_i$ are decorated by
$F_i(e_{v_i})$, $1 \leq i \leq k$, the remaining vertices are decorated
by $\beta(e_v)$.}
\end{figure}

For bi-equivariant linear maps $\Rada F1k \in \sigmabimodmaps(E,\End_V)$,
homomorphism $\beta : \freePROP(E) \to \End_V$ and {\em
distinct\/} vertices $\Rada v1k \in \Vert(G)$, we denote by
\begin{equation}
\label{stanu_se_UL_instruktorem?}
\Gamma_{\{\beta\}}^{\{\Rada v1k\}}[\Rada F1k] \in \End_V(G)
\end{equation} 
the $\End_V$-decorated graph whose vertices $v_i$, $1 \leq i \leq k$,
are decorated by $F_i(e_{v_i})$ and the remaining vertices $v \not\in
\{\Rada v1k\}$ by $\beta(v_e)$.  See Figure~\ref{figu1}.
The \PROP\ structure of $\End_V$ 
determines  the contraction $\alpha :
\End_V(G) \to \End_V(m,n)$ along $G$~\cite[Section~8]{markl:handbook}. Applying
this contraction to~(\ref{stanu_se_UL_instruktorem?}),
we obtain a linear map
\[
\alpha(\Gamma_{\{\beta\}}^{\{\Rada v1k\}}[\Rada F1k]) \in \End_V(m,n) =
\Lin(\otexp Vn,\otexp Vm).
\]

Let us show, after these preliminaries, how the $L_\infty$-braces of
Theorem~\ref{main} can be constructed.
Assume that, as in the introduction, $\alpha : \sfP \to \End_V$ is a
$\sfP$-algebra and $\rho : (\sfM,\pa) \to (\sfP,0)$ a minimal model of
$\sfP$. Recall that $\beta$ denotes the composition $\alpha \circ \rho : \sfM
\to \End_V$.  Assume that $\sfM = \freePROP(E)$ for some $\Sigma$-bimodule
$E$. It follows from definition~(\ref{Joshua_v_Praze}) and
isomorphism~(\ref{jitka2}) that
\begin{equation}
\label{piji_kavu}
\CP^*(V;V) \cong \susp  \sigmabimodmaps^{-*}(E,\End_V).
\end{equation}

For $\xi \in E(m,n)$, represent the value
$\pa(\xi) \in \freePROP(E)(m,n)$ of the differential as a sum of
$E$-decorated $(m,n)$-graphs, 
\begin{equation}
\label{porad_odkladam_ty_testy}
\pa(\xi) = \sum_{s\in S_\xi} \Gamma_s,
\end{equation}
with a finite set of summation indices $S_\xi$.
Let $F_i \in \sigmabimodmaps(E,\End_V)$ correspond, under
isomorphism~(\ref{piji_kavu}),  to a cochain $f_i \in
\CP^*(V;V)$, $1 \leq i \leq k$, and define $l_k (\Rada f1k)(\xi)  \in
\End_V(m,n)$ by
\begin{equation}
\label{nabla}
l_k (\Rada f1k)(\xi) := (-1)^{\nu(\Rada f1k)}
\sum_{s \in S_\xi}\ 
\sum_{\Rada v1k} \alpha(\Gamma_{s,{\{\beta\}}}^{\{\Rada v1k\}}[\Rada  F1k]),
\end{equation}
where $\Rada v1k$ runs over all $k$-tuples of distinct vertices of the
underlying graph $G_s$ of the $E$-decorated graph
$\Gamma_s$. The overall sign in the right
hand side, defined later
in~(\ref{dnes_jsem_byl_na_sipcich}), plays no role in this section.
The linear map $\xi \mapsto l_k (\Rada f1k)(\xi)$ determines,
by~(\ref{piji_kavu}), an element $l_k (\Rada f1k) \in \CP^*(V;V)$,
which is precisely the $k$-th $L_\infty$-bracket of Theorem~\ref{main}.
Observe that~(\ref{nabla}) makes sense also for $k=0$ when it reduces to
\[
l_0 (\xi) := 
\sum_{s \in S_\xi}\ \alpha(\Gamma_{s,{\{\beta\}}})
\]
where $\Gamma_{s,{\{\beta\}}}$ is the $E$-decorated graph whose
underlying graph is $G_s$ and all vertices $v$ are
decorated by $\beta(e_v)$. This clearly means that $\Gamma_{s,{\{\beta\}}} =
\beta(\Gamma_s)$, therefore $l_0 (\xi) = (\beta \circ
\pa)(\xi)$. Since $\beta \circ \pa = 0$, this implies that $l_0 = 0$.
It is equally simple to verify that $l_1$ coincides with the differential 
$\delta_\calP$ in $\CP^*(V;V)$. 

Let us explain why formula~(\ref{nabla}) indeed defines an
$L_\infty$-structure. It is not difficult to see that $l_k (\Rada
f1k)$, $k \geq 1$, have the appropriate symmetry. To understand why the
$L_\infty$-axiom recalled 
in~(\ref{jestli_jsem_se_nezblaznil}) below is satisfied,
expand the equation 
\hbox{$(\pa \circ \pa)(\xi) = 0$~into}
\begin{equation}
\label{21}
0 =
(\pa \circ \pa)(\xi) = \sum_{s \in S_\xi} \pa(\Gamma_s) = 
\sum_{s \in S_\xi} 
\sum_{v \in \Vert(\Gamma_s)}
\sum_{t \in T_{s,v}} \Gamma_{s,v,t}
\end{equation}
where $\Gamma_{s,v,t}$ is the $E$-decorated graph obtained as
follows. For $v \in \Vert(G_s)$, let 
\begin{equation}
\label{23}
\pa(e_v) = \sum_{t \in T_{s,v}}
\Gamma_{v,t},
\end{equation}
where $\Gamma_{v,t}$ are $E$-decorated graphs indexed by a finite
set~$T_{s,v}$.
The graph $\Gamma_{s,v,t}$ is then given
by replacing the $E$-decorated vertex $v$ of $\Gamma_s$ by the
$E$-decorated graph $\Gamma_{s,v}$. By~(\ref{21}),
\[
0 =
\sum_{s \in S_\xi} 
\sum_{v \in \Vert(\Gamma_s)}
\sum_{t \in T_{s,v}}
\sum_{\Rada v1k} \alpha(\Gamma_{s,v,t,{\{\beta\}}}^{\{\Rada v1k\}}[\Rada F1k])
\]
for arbitrary $\Rada F1k \in \sigmabimodmaps(E,\End_V)$ and $k \geq 1$.  
This summation can be further refined as
\begin{equation}
\label{uz_mne_z_ty_Kvety_hrabe}
0 =
\sum_{s \in S_\xi} 
\sum_{v \in \Vert(\Gamma_s)}
\sum_{t \in T_{s,v}} \sum_\sigma
\sum_{\Rada v1k}\hskip -.7em  {}^{\raisebox{.5em}{\scriptsize $\sigma$}}
\alpha(\Gamma_{s,v,t,{\{\beta\}}}^{\{\Rada v1k\}}[\Rada F1k]),
\end{equation}
where $\sigma$ runs over all $(i,k-i)$-unshuffles with $i \geq 1$
and the rightmost summation is restricted to $k$-tuples of 
distinct vertices $\Rada v1k \in \Vert(G_{s,v,t})$ such that $\rada
{v_{\sigma(1)}}{v_{\sigma(i)}}$ are vertices of the subgraph
$G_{s,v} \subset G_{s,v,t}$ and $\rada {v_{\sigma(i+1)}}{v_{\sigma(k)}}$ are
vertices of the complement of $G_{s,v}$ in~$G_{s,v,t}$.

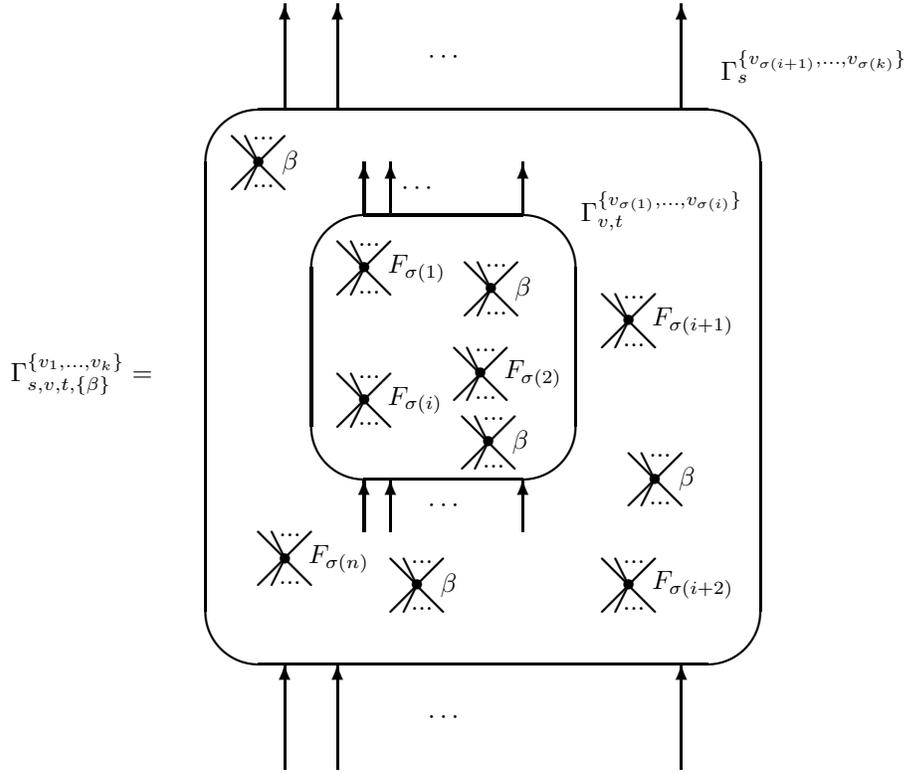
\begin{figure}[t]
\begin{center}
{
\unitlength=1.000000pt
\begin{picture}(210.00,290.00)(-10,0)
\thicklines
\put(-20,150){\makebox(0.00,0.00)[r]{$\Gamma_{s,v,t,{\{\beta\}}}^{\{\Rada
              v1k\}} =$}}
\put(107,124){\makebox(0.00,0.00){$\hnida{\beta}$}}
\put(60.00,140.00){\makebox(0.00,0.00){$\hnida{F_{\sigma(i)}}$}}
\put(195.00,260.00){%
             \makebox(0.00,0.00)[lb]%
             {$\Gamma_s^{\{\rada {v_{\sigma(i+1)}}{v_{\sigma(k)}}\}}$}}
\put(90.00,100.00){\makebox(0.00,0.00){$\cdots$}}
\put(120.00,90.00){\vector(0,1){20.00}}
\put(70.00,90.00){\vector(0,1){20.00}}
\put(60.00,90.00){\vector(0,1){20.00}}
\put(120.00,210.00){\vector(0,1){20.00}}
\put(90.00,160.00){\oval(100.00,100.00)}
\put(108,182){\makebox(0.00,0.00){$\hnida{\beta}$}}
\put(20.00,230.00){\makebox(0.00,0.00){$\hnida{\beta}$}}
\put(80.00,70.00){\makebox(0.00,0.00){$\hnida{\beta}$}}
\put(170.00,110.00){\makebox(0.00,0.00){$\hnida{\beta}$}}
\put(30.00,80.00){\makebox(0.00,0.00){$\hnida{F_{\sigma(n)}}$}}
\put(160.00,70.00){\makebox(0.00,0.00){$\hnida{F_{\sigma(i+2)}}$}}
\put(142,204){\makebox(0.00,0.00)[lb]{$\Gamma_{v,t}^{\{\rada
                                {v_{\sigma(1)}}{v_{\sigma(i)}}\}}$}}
\put(160.00,170.00){\makebox(0.00,0.00){$\hnida{F_{\sigma(i+1)}}$}}
\put(104,150.00){\makebox(0.00,0.00){$\hnida{F_{\sigma(2)}}$}}
\put(60.00,190.00){\makebox(0.00,0.00){$\hnida{F_{\sigma(1)}}$}}
\put(80.00,220.00){\makebox(0.00,0.00){$\cdots$}}
\put(90.00,20.00){\makebox(0.00,0.00){$\cdots$}}
\put(90.00,270.00){\makebox(0.00,0.00){$\cdots$}}
\put(70.00,210.00){\vector(0,1){20.00}}
\put(60.00,210.00){\vector(0,1){20.00}}
\put(180.00,0.00){\vector(0,1){40.00}}
\put(50.00,0.00){\vector(0,1){40.00}}
\put(30.00,0.00){\vector(0,1){40.00}}
\put(180.00,250.00){\vector(0,1){40.00}}
\put(50.00,250.00){\vector(0,1){40.00}}
\put(30.00,250.00){\vector(0,1){40.00}}
\put(105.00,145.00){\oval(210.00,210.00)}
\end{picture}}
\end{center}
\caption{\label{zatim_pocitac_drandi}%
The graph $\Gamma_{s,v,t,{\{\beta\}}}^{\{\Rada v1k\}}$ obtained by
replacing the vertex $v$ of $\Gamma_s^{\{\rada
{v_{\sigma(i+1)}}{v_{\sigma(k)}}\}}$ by $\Gamma_{v,t}^{\{\rada
{v_{\sigma(1)}}{v_{\sigma(i)}}\}}$.
}
\end{figure}

It is obvious that for such $\sigma$ and $\Rada v1k$, the graph
$\Gamma_{s,v,t,{\{\beta\}}}^{\{\Rada v1k\}}$ is obtained from the
$\End_V$-decorated graph $\Gamma_s^{\{\rada
{v_{\sigma(i+1)}}{v_{\sigma(k)}}\}}$ by replacing the vertex $v$ by the
$\End_V$-decorated graph $\Gamma_{v,t}^{\{\rada
{v_{\sigma(1)}}{v_{\sigma(i)}}\}}$, see
Figure~\ref{zatim_pocitac_drandi}.
Therefore one can reinterpret the right hand side
of~(\ref{uz_mne_z_ty_Kvety_hrabe}) as
\begin{equation}
\label{jestli_jsem_se_nezblaznil}
0 = \sum_{i+j = k+1} \sum_\sigma \eta(\sigma) (-1)^{i(j-1)} \cdot
l_j(l_i(\Rada f{\sigma(1)}{\sigma(i)}),\Rada f{\sigma(i+1)}{\sigma(k)})
\end{equation}
with $\eta(\sigma) := {\it sgn\/} (\sigma) \cdot \epsilon(\sigma)$,
which is the axiom of $L_\infty$-algebras.

We are sure that the reader has already realized at which points we
were not precise. First, we did not say what is a decoration
of a graph. Second, our formulas~(\ref{porad_odkladam_ty_testy})
and~(\ref{23}) for the differential assumed choices of representatives
of decorated graphs, and a rigorous proof
of~(\ref{jestli_jsem_se_nezblaznil}) would require
assumptions about the compatibility of these choices. We also
ignored signs. Namely the compatibility assumption would
make a rigorous version of the above arguments very complicated.


\section{Pre-Lie structures on spaces of derivations}
\label{2}

In Theorem~\ref{Ji} of this section we prove that, for a free \PROP\
$\sfM = \freePROP(E)$ and for an arbitrary \PROP\ $\calE$, the space $\DER$
of derivations of $\sfM$ with values in the coproduct $\sfM \copr
\calE$ admits a natural pre-Lie algebra structure.  Observe that if
$\calE$ is the trivial \PROP, then $\DER = \Der(\sfM)$ and
Theorem~\ref{Ji} is an analog of the classical theorem about the
existence of a pre-Lie structure on the space of (co)derivations of a (co)free
algebra~\cite[Section~II.3.9]{markl-shnider-stasheff:book}.
We will also study how this pre-Lie structure behaves with respect to
some natural maps induced by a \PROP\ homomorphism $\beta :\sfM \to \calE$
(Lemma~\ref{nahravam_Bacha}).

Recall that if \PROP{s} $\sfP_1$ and $\sfP_2$ are represented
as quotients of free \PROP{s}, $\sfP_s = \freePROP(X_s)/(R_s)$,
$s=1,2$, their coproduct $\sfP_1 \copr \sfP_2$ is the quotient
$\freePROP(X_1,X_2)/(R_1,R_2)$, where $(R_1,R_2)$ denotes the
\PROP{ic} ideal generated by $R_1 \cup R_2$.  The following technical
proposition will be useful in the sequel.

\begin{proposition}
\label{JITb}
Given \PROP{s} $\sfP_1$, $\sfP_2$ and a $\sfP_1 * \sfP_2$-module $\modU$,
there is a canonical isomorphism
\begin{equation}
\label{jitka1}
\Der(\sfP_1 \copr \sfP_2,\modU) \cong \Der(\sfP_1,\modU)
\oplus \Der(\sfP_2,\modU)
\end{equation}
which sends $\theta \in \Der(\sfP_1 \copr \sfP_2,\modU)$ into the
direct sum $\theta|_{\sfP_1} \oplus \theta|_{\sfP_2}$ of
restrictions. In the right hand side of~(\ref{jitka1}), the
$\sfP_i$-module structure on $\modU$ is induced from the
$\sfP_1 \copr \sfP_2$-structure by
the inclusion $\sfP_i \hookrightarrow \sfP_1 \copr \sfP_2$,
$i=1,2$.
\end{proposition}

\noindent 
The {\it proof\/} follows from the
representation~(\ref{jitka3}) of derivations as homomorphisms and the
universal property of coproducts.%
\qed

The last thing we need to observe before coming to the main point of
this section is that, given a homomorphism $\omega : \modU' \to
\modU''$ of $\sfP$-modules and a derivation $\theta \in
\Der(\sfP,\modU')$, the composition $\omega \circ \theta$ of
$\Sigma$-bimodule maps is a derivation in
$\Der(\sfP,\modU'')$. The correspondence $\theta \mapsto
\omega \circ \theta$ therefore induces the `standard' map
\begin{equation}
\label{to_jsem_zvedav_na_Jitku} \omega_* :\Der(\sfP,\modU') \to
\Der(\sfP,\modU'').
\end{equation}

Suppose that $\sfM$ and $\calE$ are \PROP{s}. The central object
of this section will be the graded vector space
$\Der(\sfM,\sfM \copr \calE)$, where the coproduct $\sfM \copr \calE$ is
considered as an $\sfM$-module with the structure induced by
the canonical inclusion
$i_\sfM : \sfM \hookrightarrow \sfM \copr \calE$. We need to introduce,
for the purposes of the next section, three maps $A$, $B$ and $C$ that
relate $\Der(\sfM,\sfM \copr \calE)$ with other spaces of derivations.

Suppose that $\calE$ is equipped with a \PROP\ homomorphism $\beta :
\sfM \to \calE$.  It then makes sense to consider
$\Der(\sfM,\calE_\beta)$, where $\calE_\beta$ denotes $\calE$ with the
$\sfM$-module structure given by the homomorphism $\beta$.  The map
$\beta$ also induces a \PROP\ homomorphism $\whb : \sfM \copr \calE \to
\calE$ by $\whb |_\sfM := \beta$ and $\whb |_\calE := \iden_\calE$. By
the definition of $\calE_\beta$, $\whb$ can be considered as 
a map of $\sfM$-modules which
in turn induces the standard map~(\ref{to_jsem_zvedav_na_Jitku})
\begin{equation}
\label{A}
C:= \C : \Der(\sfM,\sfM \copr \calE) \to \Der(\sfM,\calE_\beta).
\end{equation}
Similarly, the canonical inclusion $i_\sfM : \sfM \hookrightarrow
\sfM \copr \calE$ induces an inclusion of vector spaces
\begin{equation}
\label{B}
B := \B : \Der(\sfM) \hookrightarrow
\Der(\sfM,\sfM \copr \calE).
\end{equation}
From this moment on, 
we suppose that $\sfM$ is the {\em free \PROP\ $\sfM = \freePROP(E)$
generated by a $\Sigma$-bimodule $E$\/}.
Let $i_\calE : \calE \hookrightarrow \sfM \copr \calE$ be the
canonical inclusion and denote by $A$ the composition
\begin{equation}
\label{C}
A :
\Der(\sfM,\calE_\beta) \cong \sigmabimodmaps(E,\calE)
\stackrel{{i_\calE}_*}{\longrightarrow}
\sigmabimodmaps(E,\sfM \copr \calE)
\cong \Der(\sfM,\sfM \copr \calE),
\end{equation}
with the isomorphisms given by Proposition~\ref{JITa}.
The three maps introduced above can be organized into the diagram
\begin{equation}
\label{leze_na_mne_neco}
\Der(\sfM) \stackrel{B}{\hookrightarrow} \hskip 6em
\dvojicezob{\Der(\sfM,\sfM \copr \calE)}{\Der(\sfM,\calE_\beta)}
{C}{A}  \hskip 4.3em. \hskip 1em \rule{0pt}{1.8em}
\end{equation}
\rule{0pt}{1.8em}We will use the inclusion $B$ to identify
$\Der(\sfM)$ with a subspace of $\Der(\sfM,\sfM \copr \calE)$.
The following simple lemma will be useful.

\begin{lemma}
\label{pred_odjezdem_do_Srni}
The linear maps defined in~(\ref{A})--(\ref{C}) above  satisfy
\begin{eqnarray}
\label{jedna}
\CA &=& \iden : \Der(\sfM,\calE_\beta)
\to\Der(\sfM,\calE_\beta)\ \mbox { and }
\\
\label{dve}
\CBB &=& \beta_*: \Der(\sfM) \to \Der(\sfM,\calE_\beta).
\end{eqnarray}
\end{lemma}

\begin{proof}
By definition, for $F \in \Der(\sfM,\calE_\beta)$,
$\CA(F)|_E = \whb \circ i_\calE \circ F|_E = F|_E$, because $\whb
\circ i_\calE = \iden_{\calE}$ by the definition of $\whb$.  Since
each derivation in $\Der(\sfM,\calE_\beta)$ is, by
Proposition~\ref{JITa}, determined by its restriction to the space
of generators, this proves~(\ref{jedna}).  Similarly, for $\Phi \in
\Der(\sfM)$, $\CBB(\Phi)|_E = \whb \circ i_\sfM \circ \Phi|_E = \beta \circ
\Phi|_E$, again by the definition of $\whb$, which
proves~(\ref{dve}).
\end{proof}

Before we formulate the next statement, we observe that
Proposition~\ref{JITb} implies
\begin{equation}
\label{Jitk}
\Der(\sfM,\sfM \copr \calE) \cong \{\widetilde \theta \in \Der(\sfM \copr
\calE);\ \w\theta(\calE) = 0\}.
\end{equation}
In words, each derivation $\theta \in \Der(\sfM,\sfM \copr \calE)$
can be uniquely extended into a derivation
$\widetilde \theta \in \Der(\sfM \copr
\calE)$ characterized by $\widetilde\theta|_\calE = 0$.

\begin{lemma}
Let  $\phi, \psi \in \Der(\sfM,\sfM \copr \calE)$ be two derivations.
Then the commutator of the
composition
\begin{equation}
\label{Jitka_zvala}
[\phi,\psi] := \w \phi \circ \psi -
(-1)^{|\phi||\psi|} \cdot \w \psi \circ\phi
\end{equation}
is again a derivation, and the assignment $\phi,\psi \mapsto
[\phi,\psi]$ makes $\Der(\sfM,\sfM \copr \calE)$ a graded Lie algebra.
\end{lemma}

\begin{proof}
It is straightforward to check that $[\phi,\psi]$ defined
in~(\ref{Jitka_zvala}) has the derivation property with respect to
both the vertical and horizontal compositions. The rest of the lemma
is obvious.
\end{proof}

Let us look more closely at the isomorphism
\begin{equation}
\label{zitra_prohlidka} \Der(\sfM,\sfM \copr \calE) \cong
\sigmabimodmaps(E,\sfM \copr \calE)
\end{equation}
that follows from Proposition~\ref{JITa}.
It sends $\theta \in \Der(\sfM,\sfM \copr \calE)$ into the
restriction $\theta|_E \in
\sigmabimodmaps(E,\sfM \copr \calE)$. In the opposite
direction, to each $u \in \sigmabimodmaps(E,\sfM \copr
\calE)$ there exists a unique extension $\E(u) \in
\Der(\sfM,\sfM \copr \calE)$ characterized by
\hbox{$\E(u)|_E = u$}.

\begin{theorem}
\label{Ji}
The graded vector space $\Der(\sfM,\sfM \copr
\calE)$ is a natural graded pre-Lie algebra, with the structure
operation $\diamond : \Der(\sfM,\sfM \copr \calE)\ot \Der(\sfM,\sfM \copr \calE)
\to  \Der(\sfM,\sfM \copr \calE)$ given by
\begin{equation}
\label{Jitka!}
\theta \di \phi := (-1)^{|\theta| |\phi|}\cdot \E(\w\phi \circ
\theta|_E),
\end{equation}
where $\E(\w\phi \circ \theta|_E) \in \Der(\sfM,\sfM \copr \calE)$
is the extension of the composition
\[
\w\phi \circ \theta|_E: E
\stackrel{\theta|_E}{\longrightarrow} \sfM \copr \calE
\stackrel{\w\phi}{\longrightarrow} \sfM \copr \calE \in
\sigmabimodmaps(E,\sfM \copr \calE).
\]
\end{theorem}

\noindent 
{\it Proof\/}
of Theorem~\ref{Ji} is straightforward, but since this theorem is
a central technical tool of this section, 
we give it here. For the ease of reading, we omit in this proof
the $\w{\hskip .8em}$
denoting the extension of derivations of $\Der(\sfM,\sfM \copr \calE)$ into
derivations of $\Der(\sfM \copr \calE)$. By
definition~\cite[Section~2]{gerstenhaber:AM63}, $\di$~is a
pre-Lie product if the associator
\[
A(\theta, \phi , \psi) :=
(\theta \di \phi) \di \psi - \theta \di (\phi \di \psi)
\]
is (graded) symmetric in $\phi$ and $\psi$. By~(\ref{Jitka!}),
this associator can be written as
\begin{eqnarray}
\label{Jit}
A(\theta, \phi , \psi) &=& (-1)^{\epsilon}
\left\{
\E(\psi \circ \E(\phi \circ \theta|_E)|_E) -
\E(\E (\psi \circ \phi|_E) \circ \theta|_E)
\right\}
\\
&=& (-1)^{\epsilon} \nonumber
\left\{
\E(\psi \circ \phi \circ \theta|_E) -
\E(\E (\psi \circ \phi|_E) \circ \theta|_E)
\right\},
\end{eqnarray}
where $\epsilon := |\phi||\psi| +|\phi||\theta| +|\psi||\theta|.$
Since $A(\theta,\phi,\psi)$ is a derivation belonging to
$\Der(\sfM,\sfM \copr \calE)$, it is determined by is restriction to $E$. By~(\ref{Jit}),
clearly
\begin{equation}
\label{Jitula}
A(\theta,\phi,\psi)|_E = (-1)^{\epsilon}
\left\{
\psi \circ \phi \circ \theta|_E - \E(\psi \circ \phi|_E) \circ \theta|_E
\right\}.
\end{equation}
The antisymmetry of $A(\theta,\phi,\psi)$ in $\phi$ and $\psi$ is then
equivalent to the antisymmetry of the restrictions to $E$,
\[
A(\theta,\phi,\psi)|_E = (-1)^{|\phi||\psi|}
A(\theta,\psi,\phi)|_E
\]
which is, by~(\ref{Jitula}), the same as
\[
\psi \circ \phi \circ \theta|_E - \E(\psi \circ \phi|_E) \circ
\theta|_E -
(-1)^{|\phi||\psi|} \left\{
\phi \circ \psi \circ \theta|_E - \E(\phi \circ
\psi|_E) \circ \theta|_E
\right\} = 0,
\]
where we, of course, omitted the overall factor $(-1)^\epsilon$.
Using the bracket~(\ref{Jitka_zvala}) and moving $\theta|_E$ to the right,
the above display can be rewritten as
\[
\left\{ [\psi,\phi] - \E(\psi \circ \phi|_E) + (-1)^{|\phi||\psi|}\E(\phi \circ
\psi|_E) \right\} \circ\theta|_E = 0
\]
so it is enough to prove that
\[
[\psi,\phi] - \E(\psi \circ \phi|_E)+
(-1)^{|\phi||\psi|}\E(\phi \circ \psi|_E)
 = 0.
\]
Since the left hand side is an element of $\Der(\sfM,\sfM \copr \calE)$,
it suffices to prove that it vanishes when restricted to
generators, that is
\[
[\psi,\phi]|_E - \psi \circ \phi|_E + (-1)^{|\phi||\psi|} \phi \circ
\psi|_E = 0,
\]
which immediately follows from the definition~(\ref{Jitka_zvala}) of the
bracket.%
\qed

The last statement in this section relates the $\di$-product of
Theorem~\ref{Ji} with the
maps $A$ and $B$.

\begin{lemma}
\label{nahravam_Bacha}
Let $A$ and $C$ be the maps defined in~(\ref{A}) and~(\ref{C}). Then
for each $\theta, \phi \in \DER$,
\begin{equation}
\label{pisu_v_Ratajich}
AC (\theta) \di \phi = 0\ \mbox{ and }\ C( \theta \di AC(\phi)) =
C( \theta \di \phi).
\end{equation}
\end{lemma}

\begin{proof}
By Proposition~\ref{JITa}, it suffices to
prove the restrictions of the above equalities onto the space $E$ of
generators of $\sfM = \freePROP(E)$. By definition,
\[
(AC(\theta) \di \phi)|_E 
= (-1)^{|\phi||\theta|} \cdot \w \phi \circ AC(\theta)|_E 
= (-1)^{|\phi||\theta|}\cdot\w \phi \circ i_\calE \circ \whb
\circ \theta|_E = 0,
\]
because $\w \phi \circ i_\calE = \w \phi|_\calE = 0$. This proves the
left equation of~(\ref{pisu_v_Ratajich}).
Similarly, by definition 
\[
C(\theta \di AC(\phi))|_E
= (-1)^{|\phi||\theta|} \cdot
\whb \circ \widetilde{AC(\phi)} \circ \theta|_E.
\]
Let us prove that $\whb \circ \w{AC(\phi)} =
\whb \circ \w \phi$ in $\Der(\sfM
\copr \calE)$. By Proposition~\ref{JITb} this means to verify that
\[
\whb \circ \w{AC(\phi)}|_\calE = \whb \circ \w \phi|_\calE\
\mbox { and }\ \whb \circ \w{AC(\phi)}|_\sfM = 
\whb \circ \w \phi|_\sfM.
\] 
The first equation is obvious because $\w{ \AC(\phi)}|_\calE = 0 =
\w {\phi}|_\calE$ by the definition of the tilde-extension. The second
equality is established by
\[
\whb \circ \w{\AC(\phi)}|_\sfM =
\whb \circ \AC(\phi) = \CAC(\phi) = C(\phi) = \whb \circ \w \phi|_\sfM,
\] 
where we used the definition~(\ref{A}) of the map $C$ and the equality
$\CA= \iden$ proved in Lemma~\ref{pred_odjezdem_do_Srni}. This
finishes the proof of the right equation of~(\ref{pisu_v_Ratajich}).
\end{proof}

\section{Braces}
\label{2bis}

In this section we prove Theorem~\ref{main} and
Proposition~\ref{main1} of the introduction and formulate
some technical statements which will guarantee the convergence of the
master equation.

According to~\cite{guin-oudom,lm:sb}, the pre-Lie algebra product
$\di$ on $\DER$ whose existence we proved in 
Theorem~\ref{Ji}
generates unique {\em symmetric braces\/}. This means
that for each $\theta,\phi_1,\ldots, \phi_n \in \Der(\sfM,\SFM \copr
\calE)$, $n \geq 1$, there exists a `brace'
$\theta\langle\phi_1,\ldots, \phi_n\rangle \in \Der(\sfM,\SFM \copr
\calE)$ such that 
\begin{equation}
\label{dnes_to_zkusim_nainstalovat}
\theta \langle \phi_1 \rangle = \theta \di \phi_1.
\end{equation}
These braces satisfy the axioms recalled in the
Appendix~A on page~\pageref{apA}, 
where we also indicate how these braces are generated by
$\di$. The following statement generalizes Lemma~\ref{nahravam_Bacha}.

\begin{lemma}
\label{opravu_se_na_nem_neda_pracovat}
Let $A$ and $C$ be the maps defined in~(\ref{A}) and~(\ref{C}). Then
for each $\theta,\Rada \phi 1n \in Der(\sfM,\SFM \copr \calE)$, $n
\geq 1$,
\begin{equation}
\label{Monteverdi}
AC\theta \langle \Rada \phi1n \rangle = 0
\end{equation}
and
\begin{equation}
\label{Claudio}
C\lv \theta \langle \Rada \phi1n \rangle\rv =
C\lv\theta \langle \Rada {AC\phi}1n \rangle\rv.
\end{equation}
\end{lemma}

Before we prove the lemma we notice 
that~(\ref{Monteverdi}) for $n=1$ (with $\phi = \phi_1$) says that
\begin{equation}
\label{Jitulka}
AC\theta \langle \phi \rangle = 0
\end{equation}
which is, by~(\ref{dnes_to_zkusim_nainstalovat}), the same as 
$AC (\theta) \di \phi = 0$, which we recognize as the first
equality in~(\ref{pisu_v_Ratajich}). 
Similarly,~(\ref{Claudio}) for $n=1$ means that
\begin{equation}
\label{Pincakova}
C\lv \theta \langle \phi \rangle \rv = C\lv \theta \langle AC \phi \rangle\rv
\end{equation}
which is, again by~(\ref{dnes_to_zkusim_nainstalovat}), the same as
$C( \theta \di AC(\phi)) =
C( \theta \di \phi)$, the second equality in~(\ref{pisu_v_Ratajich}).
Therefore Lemma~\ref{opravu_se_na_nem_neda_pracovat} 
indeed generalizes Lemma~\ref{nahravam_Bacha}.

\begin{proof}[Proof of Lemma~\ref{opravu_se_na_nem_neda_pracovat}]
We prove~(\ref{Monteverdi}) by induction. For $n=1$ it
is~(\ref{Jitulka}). Assume we have already
proved~(\ref{Monteverdi}) for all $1\leq n < N$ and prove it for $n=N$. 
By axiom~(\ref{axiom}) of symmetric braces, 
\[
AC\theta \langle \sqRada \phi1N \rangle \hskip - 0.5pt =\hskip - 0.5pt
AC\theta \langle \phi_1 \rangle\langle \sqRada \phi 2N \rangle
- 
\sum \epsilon   \cdot 
AC\theta \langle \phi_1 \langle 
\sqRada \phi {i_1}{i_a} \rangle, \sqRada \phi {j_1}{j_b} \rangle,
\]
where the sum in the right hand side runs over all unshuffles
\begin{equation}
\label{Jarka_pred_chvili_volala.}
i_1 <  \cdots < i_a, \  j_1 <  \cdots < j_b,\ a \geq 1,\ a+b = N-1,
\end{equation}
of the set $\{\rada 2N\}$. The sign $\epsilon$ is not important for
the purposes of this proof, because all terms in the right hand side
are zero, by induction. This establishes~(\ref{Monteverdi}) for $n=N$.

Equation~(\ref{Claudio}) will also be be proved by induction. For $n=1$ it
is~(\ref{Pincakova}). 
Suppose we have already established~(\ref{Claudio}) 
for all $1\leq n < N$ and prove it for $n=N$. 
By~(\ref{axiom}),  
\begin{eqnarray}
\label{Hribecek} \nonumber 
\lefteqn{C\lv \theta \langle \Rada {AC\phi}1N \rangle \rv = \hskip 2em}
\\
&&\hskip 2em - \sum \epsilon \cdot C\lv
\theta \langle AC\phi_1 \langle \Rada {AC\phi}{i_1}{i_a} \rangle,
\Rada {AC\phi}{j_1}{j_b} \rangle
\rv
\\ \nonumber 
&&\hskip 2em 
+\
C\lv
\theta \langle AC\phi_1 \rangle
\langle \Rada {AC\phi}2N \rangle\rv,
\end{eqnarray}
with the sum in the right hand side taken over the
set~(\ref{Jarka_pred_chvili_volala.}) and $\epsilon$ an appropriate sign.
The sum is zero by~(\ref{Monteverdi})
while second term equals
\[
C\lv \theta \langle AC\phi_1 \rangle
\langle \Rada {\phi}2N \rangle\rv
\]
by induction.  Using~(\ref{axiom}), we can write
\begin{eqnarray*}
\lefteqn{
C\lv \theta \langle AC\phi_1 \rangle
\langle \Rada {\phi}2N \rangle\rv = \hskip 2em}
\\
&&\hskip 2em \sum \epsilon \cdot C\lv
\theta \langle AC\phi_1 \langle \Rada {\phi}{i_1}{i_a} \rangle,
\Rada {\phi}{j_1}{j_b} \rangle
\rv
+\ C\lv \theta \langle AC\phi_1, \Rada {\phi}2N \rangle\rv,
\end{eqnarray*}
where the summation range and $\epsilon$ are the same as
in~(\ref{Hribecek}).  The sum in the right hand side is zero
by~(\ref{Monteverdi}), while the second term equals, by~(\ref{axiom}),
\begin{equation}
\label{mam_hlad}
\epsilon_N \cdot
C\lv \theta\langle \sqRada {\phi}2N \rangle \langle AC \phi_1 \rangle\rv
-
\sum_{2 \leq i \leq N} \epsilon_i \cdot C\lv \theta 
\langle \phi_2,\sqldots,\phi_i \langle AC\phi_1 \rangle,
\sqldots,\phi_N \rangle\rv, 
\end{equation}
where 
\[
\epsilon_i := (-1)^{|\phi_1|(|\phi_2| + \cdots + |\phi_i|)},\
2 \leq i \leq N.
\]
The first term in~(\ref{mam_hlad}) equals
$\epsilon_N \cdot C \lv \theta\langle \Rada {\phi}2N \rangle
\langle  \phi_1 \rangle \rv$ by~(\ref{Pincakova}), 
while the second term equals, by induction,
\begin{equation}
\label{xx}
- \sum_{2 \leq i \leq N}  \epsilon_i \cdot
C \lv \theta \langle AC\phi_2,\ldots,AC(\phi_i \langle AC\phi_1
\rangle ),
\cdots,AC\phi_n \rangle \rv.
\end{equation}
Since, by~(\ref{Pincakova}), $C \lv \phi_i \langle AC\phi_1 \rangle \rv =
C \lv \phi_i \langle \phi_1 \rangle \rv$ for $2 \leq i \leq N$,~(\ref{xx}) 
equals
\[
-\sum_{2 \leq i \leq N}
\epsilon_i \cdot 
C\lv \theta \langle AC\phi_2,\ldots,AC(\phi_i \langle \phi_1
\rangle\rv ,
\cdots,AC\phi_n \rangle)
\]
which is
\[
- \sum_{2 \leq i \leq N}  \epsilon_i \cdot
C \lv \theta \langle \phi_1,\ldots,\phi_i \langle \phi_1 \rangle,
\cdots,\phi_N \rangle\rv
\]
by induction. We therefore established that
\begin{eqnarray*}
\lefteqn{
C\lv \theta \langle \Rada {AC\phi}1N \rangle\rv = \hskip 2em}
\\
&& \hskip 2em
\epsilon_N \cdot C\lv \theta\langle \Rada {\phi}2N \rangle
\langle  \phi_1 \rangle\rv
-  \sum_{2 \leq i \leq N}  \epsilon_i \cdot
C \lv \theta \langle \phi_1,\ldots,\phi_i \langle \phi_1 \rangle,
\ldots,\phi_N \rangle\rv .
\end{eqnarray*}

By~(\ref{axiom}), the right hand side of the above display equals
$C\lv \theta \langle \Rada {\phi}1N \rangle\rv$, which
establishes~(\ref{Claudio}) for $n=N$.
\end{proof}

In the following important definition, $C$ is the map introduced
in~(\ref{A}). Recall also that we use the inclusion $B$ defined
in~(\ref{B}) to identify $\Der(\sfM)$ with a subspace of $\DER$.

\begin{definition}
\label{zase_mne_jakobi_boli_hlava}
\hskip .8em For $\Phi \in \Der(\SFM)$ and $F_1,\ldots,F_n \in
\Der(\SFM,\calE_\beta)$, define the derivation 
$\Phi [ F_1,\ldots,F_n ] \in
\Der(\sfM,\calE_\beta)$ by
\begin{equation}
\label{brace}
\Phi [ F_1,\ldots,F_n ]  :=
C\left(\Phi \langle AF_1,\ldots,AF_n \rangle\right),
\end{equation}
where $\Phi \langle AF_1,\ldots,AF_n \rangle$ in the r.h.s.\
is the symmetric brace in $\Der(\sfM,\SFM \copr \calE)$.
\end{definition}

The following proposition shows that $\Der(\SFM,\calE_\beta)$ 
behaves as a left
module over the symmetric brace algebra $\Der(\SFM)$.

\begin{proposition}
\label{E_napsala}
The brace $\Phi[\Rada F1n]$ is graded symmetric in $F_1,\ldots,F_n \in
\Der(\SFM,\calE)$. Moreover, for each $\Phi_1,\ldots,\Phi_m \in \Der(\SFM)$
with $\beta \circ \Phi_1 = \cdots =\beta \circ \Phi_m = 0$,
\begin{eqnarray}
\label{7}
\lefteqn{
\Phi\langle \Rada \Phi 1m \rangle [\Rada F1n]
= \hskip 2em}
\\
&&\hskip 2em
\nonumber
\sum \epsilon \cdot
\Phi [\Phi_1 [F_{i_1^1},\ldots,F_{i_{t_1}^1}],\ldots,
\Phi_n [F_{i_1^m},\ldots,F_{i_{t_m}^m}],
F_{i_1^{m+1}},\ldots,F_{i_{t_{m+1}}^{m+1}}],
\end{eqnarray}
where the sum is taken over all unshuffle decompositions
\[
i_1^1<\cdots<i_{t_1}^1,\ldots,i_1^{m+1}<\cdots<i_{t_{m+1}}^{m+1},
\ \Rada t1{m} \geq 1,\ t_{m+1} \geq 0,
\]
of $\{1,\ldots,n\}$ and where $\epsilon$ is the Koszul sign of the
corresponding permutation of $\Rada F1n$.
\end{proposition}

\begin{proof}
The graded symmetry of the brace~(\ref{brace})
immediately follows from the definition. Let us prove~(\ref{7}).
We have
\[
\Phi\langle \Rada \Phi 1m \rangle [\Rada F1n]
=
C\left(
\Phi\langle \Rada \Phi 1m \rangle \langle\rada {AF_1}{AF_n}\rangle
\right)
\]
which can be, using~(\ref{axiom}), expanded into
\[
\sum \epsilon \cdot C\left(
\Phi \langle \Phi_1  \langle AF_{i_1^1},
\sqldots,AF_{i_{t_1}^1} \rangle,\sqldots,
\Phi_n \langle AF_{i_1^m},\sqldots,AF_{i_{t_m}^m}\rangle,
AF_{i_1^{m+1}},\sqldots,AF_{i_{t_{m+1}}^{m+1}}\rangle
\right),
\]
where $\epsilon$ and the sum is the same as
in~(\ref{axiom}). 
By~(\ref{Claudio}), this equals
\begin{eqnarray*}
\sum \epsilon \cdot C\left(
\Phi \langle AC\Phi_1  \langle AF_{i_1^1},
\ldots,AF_{i_{t_1}^1} \rangle,\ldots,
AC\Phi_n \langle AF_{i_1^m},\ldots,AF_{i_{t_m}^m}\rangle,\right.
\hskip -10em&&
\\
&&\left. 
\ACA F_{i_1^{m+1}},\ldots,\ACA F_{i_{t_{m+1}}^{m+1}}\rangle
\right)
\end{eqnarray*}
which, by definition~(\ref{brace}) of the braces, can be rewritten as
\begin{eqnarray}
\label{prohlidkou_jsem_prosel!}
\sum \epsilon \cdot
\Phi [ C\Phi_1  \langle AF_{i_1^1},\ldots,AF_{i_{t_1}^1} \rangle,\ldots,
C\Phi_n \langle AF_{i_1^m},\ldots,AF_{i_{t_m}^m}\rangle, 
\hskip -8em&&
\\
&& \nonumber 
\CA F_{i_1^{m+1}},\ldots,\CA F_{i_{t_{m+1}}^{m+1}}].
\end{eqnarray}

At this point we need to 
observe that, by~(\ref{dve}), $\beta \circ \Phi_j = C \Phi_j$ (recall
that we identified $\Phi_j$ with its image $B\Phi_j$). Therefore 
the assumption $\beta \circ \Phi_j = 0$ 
for $1 \leq j \leq m$ implies that we may assume, in the
sum~(\ref{prohlidkou_jsem_prosel!}), that all $t_j \geq 1$, because,
if $t_j=0$,
\[
C\Phi_j  \langle AF_{i_1^j},\ldots,AF_{i_{t_j}^j} \rangle = 
C \Phi_j \langle \hskip .5em \rangle =
C\Phi_j = \beta \circ \Phi_j = 0.
\] 
Since $\CA = \iden$ by~(\ref{jedna}), the term 
in~(\ref{prohlidkou_jsem_prosel!}) equals 
\[
\sum \epsilon \cdot
\Phi [ C\Phi_1  \langle AF_{i_1^1},\ldots,AF_{i_{t_1}^1} \rangle,\ldots,
C\Phi_n \langle AF_{i_1^m},\ldots,AF_{i_{t_m}^m}\rangle,
F_{i_1^{m+1}},\ldots,F_{i_{t_{m+1}}^{m+1}}],
\]
with the same summation as in~(\ref{7}).
By the definition~(\ref{brace}) of the braces, this is precisely
the right hand side of~(\ref{7}).
\end{proof}

As usual,  $\susp W$ (resp.~$\desusp
W$) denotes the suspension (resp.~desuspension) of a graded vector
space $W$. We use the same symbols to denote also the corresponding
maps $\uparrow : \hskip .2em \desusp W \to W$ and $\downarrow : W \to
\hskip .2em \desusp W$. In the following theorem,
$\Rada f1n$ will be elements of $\susp \Der(\SFM,\calE_\beta)$
and
\begin{equation}
\label{dnes_jsem_byl_na_sipcich}
\nu(\Rada f1n) :=  (n-1)|f_1| + (n-2)|f_2| + \cdots + |f_{n-1}|.
\end{equation}

\begin{theorem}
\label{E}
Let $\pa \in \Der(\SFM)$ be a degree $-1$ derivation such that $\pa^2=0$
and $\beta \circ\pa = 0$. Then the formula
\begin{equation}
\label{zpatky_v_Praze}
l_n (\Rada f1n) := (-1)^{\nu(\Rada f1n)}
\cdot  \susp  \pa [ \Rada {\desusp f}1n ]
\end{equation}
defines on the suspension $\susp \Der(\SFM,\calE_\beta)^{-*}$ a structure of an
$L_\infty$-algebra{\rm~\cite{lada-markl:CommAlg95}}.
\end{theorem}

\begin{proof}
Observe first that $\pa^2=0$ is,
by~(\ref{dnes_to_zkusim_nainstalovat}), equivalent to $\pa\langle \pa
\rangle = 0$. Expanding
\[
0 = (-1)^{\nu(\Rada f1n)} \cdot
\pa \langle \pa \rangle [ \Rada {\desusp f}1n ]
\]
using~(\ref{7}) we obtain
\[
0= \sum_{i+j = n+1} \sum_\sigma  \epsilon(\sigma)
(-1)^{\nu(\Rada f1n)} \cdot
\pa[\pa[\Rada {\desusp \hskip -.2em f}{\sigma(1)}{\sigma(i)}],\Rada
  {\desusp \hskip -.2em f}{\sigma(i+1)}{\sigma(n)}]
\]
with $\sigma$ running over all $(i,n-i)$-unshuffles with $i
\geq 1$ and $\epsilon(\sigma)$ the Koszul sign of the permutation
\[
\Rada f1n \mapsto \Rada f{\sigma(1)}{\sigma(n)}.
\]
Substituting for $l_i$ and $l_j$ from~(\ref{zpatky_v_Praze}) gives
\begin{equation}
\label{dnes_mne_strasne_bolela_zada}
0 = \sum_{i+j = n+1} \sum_\sigma \eta(\sigma) (-1)^{i(j-1)} \cdot
l_j(l_i(\Rada f{\sigma(1)}{\sigma(i)}),\Rada f{\sigma(i+1)}{\sigma(n)}),
\end{equation}
where
$\eta(\sigma)$ is as~(\ref{jestli_jsem_se_nezblaznil}).
We recognize~(\ref{dnes_mne_strasne_bolela_zada}) as the defining
equation for
$L_\infty$-algebras, see~\cite[Definition~2.1]{lada-markl:CommAlg95}.
\end{proof}

\noindent
{\it Proof of Theorem~\ref{main}\/} easily follows from Theorem~\ref{E}
applied to the situation when $\SFM$ is a minimal model $(\sfM,\pa)$
of the \PROP\ $\sfP$, $\calE = \End_V$ and $\beta = \alpha \circ \rho$
as in Section~\ref{intro}.  The condition $\pa \circ \beta =0$ is
implied by the minimality of the differential~$\pa$.%
\qed

The following proposition compares the braces defined above with the
constructions of Section~\ref{1bis}.

\begin{proposition}
\label{maminka_si_zlomila_nohu}
The $L_\infty$-structure of
Theorem~\ref{main} has the form~(\ref{nabla}) of 
Section~\ref{1bis}.
\end{proposition}

\begin{proof}
Let, in this proof, an $(E,\calE)$-decorated graph means a graph
with vertices decorated either by $E$ or by $\calE$.
We use the convention that $\Upsilon$
with a subscript will denote an $(E,\calE)$-decorated graph, and
$Y$ with the same subscript the underlying un-decorated graph. For
such $\Upsilon$, let $\vert_E(\Upsilon)$ be the set of
$E$-decorated vertices of~$\Upsilon$.

We start the proof by giving an explicit formula for the operation $\di$
introduced in Theorem~\ref{Ji}. Let $\theta, \psi \in \Der(\sfM,\sfM
\copr \calE)$ as in~(\ref{Jitka!}). 
It follows 
from~(\ref{jsem_zvedav_jestli_se_ta_druha_Jitka_ozve}) and the
definition of the free product that, 
for $\xi \in E(m,n)$, $\theta(\xi) \in \sfM \copr \calE$ is the summation
\[
\theta(\xi) = \sum_{s \in R_\xi}  \Upsilon_s
\]
of $(E,\calE)$-decorated $(m,n)$-graphs over a finite indexing set $R_\xi$.
The derivation property of $\phi$ implies that
\begin{equation}
\label{Kdy_se_nit_zivota_pretrhne?}
(\theta \di \phi)(\xi) = \sum_{s \in R_\xi} \sum_{v}
(-1)^{|\theta||\phi|} \cdot 
\Upsilon_s^{\{v\}}[\phi],
\end{equation}
where the second summation runs over $\Vert_E(\Upsilon_s)$ 
and $\Upsilon_s^{\{v\}}[\phi]$ denotes $ \Upsilon_s$ with the
decoration $e_v \in E$ of $v \in \Vert_E(\Upsilon_s)$ 
changed to $\phi(e)$. 
Let us prove inductively that the symmetric
brace induced by $\di$ satisfies
\begin{equation}
\label{zatim_nevim}
\theta \langle \rada {\phi_1}{\phi_n}\rangle =  
\sum_{s \in R_\xi} \sum_{v_1,\ldots,v_n} (-1)^{\epsilon_n} \cdot
\Upsilon_s^{\{{v_1,\ldots,v_n\}}}[\rada{\phi_1}{\phi_n}],
\end{equation}
where  $\Rada v1n$ runs over distinct elements of
$\vert_E(\Upsilon_s)$,
$\Upsilon_s^{\{{v_1,\ldots,v_n\}}}[\rada{\phi_1}{\phi_n}]$ denotes
$\Upsilon_s$ with the decoration $e_{v_i}$ of $v_i$ changed to
$\phi(v_i)$, and
\[
\epsilon_i := \zn i,
\] 
for $1 \leq i \leq n$. Since,
by~(\ref{dnes_to_zkusim_nainstalovat}), 
$\theta \langle \phi_1 \rangle = \theta \di \phi_1$, 
(\ref{zatim_nevim}) holds for $n = 1$
by~(\ref{Kdy_se_nit_zivota_pretrhne?}). 

Before continuing, we rewrite the right hand side
of~(\ref{zatim_nevim}) into a sum of $(E,\calE)$-decorated graphs. To
this end, we introduce a notation which will be useful also later in
the proof.
Let, for $1 \leq i \leq n$,
\[
\phi_i(e_{v_i}) = \sum_{t_i \in U_{s,v_i}} \Upsilon_{t_i,v_i},
\]
where $\Upsilon_{t_i,v_i}$ are $(E,\calE)$-decorated graphs and 
$U_{s,v_i}$ a finite set. For a subset $B \subset
\{0,\ldots,n\}$, let $\Kr B$ denote the $(E,\calE)$-decorated graph
obtained from
$\Upsilon_s$ by replacing, for each $i \in B$, the  
$e_{v_i}$-decorated vertex $v_i$ by the decorated graph 
$\Upsilon_{t_i,v_i}$. With this notation, 
\[
\Upsilon_s^{\{{v_1,\ldots,v_n\}}}[\rada{\phi_1}{\phi_n}] 
=
\sum_{\Rada t1n} \Kr {\{\rada 1n\}},
\]
where $t_i$ runs over $u_{t_i,v_i}$, $1 \leq i \leq n$,
so we may rewrite the right hand side of~(\ref{zatim_nevim})~as
\begin{eqnarray}
\lefteqn{
\label{psano_v_Bonnu}
\sum_{s \in R_\xi} \sum_{v_1,\ldots,v_n}  (-1)^{\epsilon_n} \cdot
\Upsilon_s^{\{{v_1,\ldots,v_n\}}}[\phi_1,\ldots,\phi_n] 
= \hskip 10em}
\\
&& \nonumber 
\hskip 10em
\sum_{s \in R_\xi} \sum_{v_1,\ldots,v_n} 
\sum_{t_1,\ldots,t_n} 
 (-1)^{\epsilon_n} \cdot \Kr{\{\rada 1n\}}. 
\end{eqnarray}

Suppose we have proved~(\ref{zatim_nevim}) for
all $k$, $1 \leq k < n$.
Consider the equation
\begin{eqnarray}
\lefteqn{
\label{zase_jsem_nachlazeny-je_to_mozne?}
\theta \langle \rada{\phi_1}{\phi_n}\rangle
= \hskip 3em}
\\
&& \nonumber 
\hskip 3em
 \theta \langle \rada{\phi_1}{\phi_{n-1}}\rangle \langle \phi_n
\rangle 
-\sum_{1 \leq i \leq n-1} 
(-1)^{\omega}\cdot
\theta \langle \phi_1,\ldots,\phi_i \langle \phi_n \rangle, \ldots,
\phi_{n-1} \rangle,
\end{eqnarray}
which follows from axiom~(\ref{axiom}) of symmetric braces; 
in the last term
\[
\omega := |\phi_n|(|\phi_{i+1}| + \cdots + |\phi_{n-1}|).
\] 
Let us analyze the first term in the right hand side 
of~(\ref{zase_jsem_nachlazeny-je_to_mozne?}).
By the induction assumption,
\[
\theta \langle \rada{\phi_1}{\phi_{n-1}}\rangle (\xi) = 
\sum_{s \in R_\xi} \sum_{v_1,\ldots,v_{n-1}}
 (-1)^{\epsilon_{n-1}} \cdot
\Upsilon_s^{\{{v_1,\ldots,v_{n-1}\}}}[\rada{\phi_1}{\phi_{n-1}}].
\]
With the notation above,
\[
\Upsilon_s^{\{{v_1,\ldots,v_{n-1}\}}}[\rada{\phi_1}{\phi_{n-1}}]
= 
\sum_{\Rada t1{n-1}} \Kr {\{1,\ldots,n-1\}},
\]
where $t_i$ runs over $U_{s,v_i}$, $1 \leq i \leq n-1$, therefore
\begin{eqnarray}
\label{dnes_je_mne_zase_trochu_lip}
\lefteqn{
\theta \langle \rada{\phi_1}{\phi_{n-1}}\rangle \langle \phi_n
\rangle (\xi) = \hskip 5em}
\\
&& \nonumber  \hskip 5em
\sum_{s \in R_\xi} \sum_{v_1,\ldots,v_{n-1}}
\sum_{t_1,\ldots,t_{n-1}}  \sum_{v_n}  (-1)^{\epsilon_n} \cdot
\Kr {\{1,\ldots,n-1\}}^{\{v_n\}}[\phi_n],
\end{eqnarray}
where $v_n$ runs over $E$-decorated vertices of $\Kr
{\{1,\ldots,n-1\}}$. Since clearly
\begin{equation}
\label{co_doma}
\vert_E(\Kr {\{1,\ldots,n-1\}}) = 
\left(
\Vert_E(\Upsilon_s) 
\setminus   \{{v_1,\sqldots,v_{n-1}\}}
\right)
\cup \hskip -.8em
\bigcup_{1 \leq i \leq n-1}\hskip -1em \Vert_E(\Upsilon_{t_i,v_i}),
\end{equation}
the right hand side 
of~(\ref{dnes_je_mne_zase_trochu_lip}) breaks into
$n$ components,
\[
\theta \langle \rada{\phi_1}{\phi_{n-1}}\rangle \langle \phi_n
\rangle (\xi) = A_0 + A_1 + \cdots + A_{n-1},
\]
where
\[
A_0: = 
\sum_{s \in R_\xi} \sum_{v_1,\ldots,v_{n-1}}  
\sum_{t_1,\ldots,t_{n-1}} {\sum_{v_n}}^{(0)}(-1)^{\epsilon_n} \cdot
\Kr {\{1,\ldots,n-1\}}^{\{v_n\}}[\phi_n]
\]
with the superscript $(0)$ meaning that $v_n$ runs over 
$\Vert_E(\Upsilon_s) 
\setminus   \{{v_1,\ldots,v_{n-1}\}}$, and
\[
A_i: =  
\sum_{s \in R_\xi} \sum_{v_1,\ldots,v_{n-1}} 
\sum_{t_1,\ldots,t_{n-1}} {\sum_{v_n}}^{(i)}(-1)^{\epsilon_n} \cdot
\Kr {\{1,\ldots,n-1\}}^{\{v_n\}}[\phi_n],\ 1 \leq i \leq n-1,
\]
where $(i)$ means that $v_n$ runs over $\Vert_E(\Upsilon_{t_i,v_i})$.

It is immediately clear that, for $v_n \in 
\Vert_E(\Upsilon_s) 
\setminus   \{{v_1,\ldots,v_{n-1}\}}$, 
\[
\Kr{\{\rada 1{n-1}\}}^{\{v_n\}} = \sum_{t_n}  \Kr{\{\rada 1n\}}, 
\]
so
\[
A_0 =\sum_{s \in R_\xi} \sum_{v_1,\ldots,v_n}
\sum_{\Rada t1n}
(-1)^{\epsilon_n} \cdot
 \Kr{\{\rada 1n\}}, 
\]
which is, by~(\ref{psano_v_Bonnu}),
the right hand side of~(\ref{zatim_nevim}). It is equally clear that,
for $v_n \in \vert_E(\Upsilon_{t_i,v_i})$,
\[
\sum_{t_i} \sum_{v_n} 
\Kr{\{\rada 1{n-1}\}}^{\{v_n\}}
[\phi_n] =  (-1)^{\omega} \cdot
\Kr {\{1,\ldots,i-1,i+1,\ldots,n-1\}}^{\{v_i\}}[\phi_i\langle \phi_n \rangle],
\]
therefore
\begin{eqnarray*}
A_i \hskip -.2em  &=&  \hskip -.2em 
\sum_{s \in R_\xi} \sum_{v_1,\squeezedldots,v_{n-1}}
\sum_{t_1,\squeezedldots,t_{i-1}}\sum_{t_{i+1},\squeezedldots,t_{n-1}} 
 \hskip -.1em 
{\sum_{v_n}}^{(i)} \hskip -.2em 
(-1)^{\epsilon_n \hskip -.1em +\hskip -.1em\omega} 
\hskip -.2em  \cdot \hskip -.2em 
\Kr {\{1,\squeezedldots,i-1,i+1,
\squeezedldots,n-1\}}^{\{v_i\}}[\phi_i\langle \phi_n
  \rangle]
\\
&=& \hskip -.2em 
\sum_{s \in R_\xi} \sum_{v_1,\ldots,v_{n-1}}
{\sum_{v_n}}^{(i)} (-1)^{\epsilon_n +\omega}
\Upsilon_s^{\{\Rada v1{n-1}\}}[\phi_1,\ldots,\phi_i
\langle \phi_n \rangle,\ldots,\phi_{n-1}], 
\end{eqnarray*}
which equals 
$(-1)^{\omega} 
\theta \langle \phi_1,\ldots,\phi_i \langle \phi_n \rangle, \ldots,
\phi_{n-1} \rangle$, by induction. We recognize this expression as one
of the remaining terms in the right hand side
of~(\ref{zase_jsem_nachlazeny-je_to_mozne?}), taken with the minus
sign. Assembling the above results, we obtain~(\ref{zatim_nevim}).

Let $\pa(\xi) = \sum_{s\in S_\xi} \Gamma_s$ be as
in~(\ref{porad_odkladam_ty_testy}) and $F_i \in
\sigmabimodmaps(E,\End_V) \cong
\Der(M,\End_V)$, for $1 \leq i \leq n$. By~(\ref{zatim_nevim}), 
\[
\pa \langle AF_1,\ldots,AF_n \rangle (\xi) =
\sum_{s \in S_\xi} \sum_{\Rada v1n} \Gamma_s^{\{\Rada v1n\}}[\Rada F1n](\xi).
\] 
Applying, as in Definition~\ref{zase_mne_jakobi_boli_hlava}, the map 
$C$ on this identity, we get a formula for $\pa[\Rada F1n]$ which
agrees, modulo signs, to the right hand side~(\ref{nabla}). The sign
factor is induced by (de)suspensions.
\end{proof}

Proposition~\ref{maminka_si_zlomila_nohu} has several 
important implications. Let us formulate first a corollary that
implies Proposition~\ref{main1}; the notation is the same as the one
introduced in the paragraph preceding this proposition.

\begin{corollary}
\label{r}
Let $\xi \in \freePROP(E)$ be such that $\pa(\xi) \in \freePROP^{\leq
k}(E)$. Then 
\[
l_n(\Rada f1n)(\xi) = 0 \mbox { for each } n > k.
\]
In particular, if $\pa(E) \subset \freePROP^{\leq k}(E)$, then $l_n =
0$ for $n > k$.
\end{corollary}

\begin{proof}
If $\pa(\xi) \in \freePROP^{\leq
k}(E)$, then all graphs in~(\ref{porad_odkladam_ty_testy}) have $\leq k$
vertices, so the summation in~(\ref{nabla}) 
is empty for $k > 2$.%
\end{proof}

In this paper we write several formulas containing
infinite sums. Their convergence will be guaranteed by the following property
of an \li-algebra $L = (W,l_1,l_2,\ldots)$:
\begin{equation}
\label{V}
\begin{minipage}{32em}
is a direct product $W = \prod_{s \geq 1}W_s$ such that $l_k(\Rada
w1k)_s = 0$ for all $k > s$ and $\Rada w1k \in W$.
\end{minipage}
\end{equation}

In~(\ref{V}), $l_k(\Rada w1k)_s$ denotes the component of
$l_k(\Rada w1k)$ in $W_s$.  There are other
conditions that can guarantee the convergence of our formulas, as the
nilpotency~\cite[Definition~4.2]{getzler:AM:09}, but \li-algebras
considered in this paper may not be nilpotent.

\begin{proposition}
\label{conv}
The \li-structure of Proposition~\ref{main1} satisfies~(\ref{V}).
\end{proposition}

\begin{proof}
In this case
\[
W = \susp \Der(\sfM,\End_V)^{-*} \cong \prod_{m,n} \susp
\Lin_{\mbox
  {\scriptsize$\Sigma_m$-$\Sigma_n$}}(E(m,n),\End_V(m,n))^{-*}.
\] 
For each $t \geq 1$ define a
$\Sigma_m$-$\Sigma_n$-invariant subspace $U_t(m,n) \subset E(m,n)$ as
\[
U_t(m,n) := \{\xi \in E(m,n);\ \pa(\xi) \in \freePROP^{\leq t}(E)\}.
\]
Since $\Sigma_m \times\Sigma_n$ 
is finite, there are $\Sigma_m$-$\Sigma_n$-invariant subspaces
$V_s(m,n)$, $s \geq 1$, such that 
\[
U_t(m,n) = \bigoplus_{s\leq t} V_s(m,n),\ \mbox {for each $t \geq 1$}.
\]
The same reasoning as in the proof of Corollary~\ref{r} shows that the
subspaces 
\[
W_s := \prod_{m,n} \susp 
\Lin_{\mbox {\scriptsize $\Sigma_m$-$\Sigma_n$}}(V_s(m,n),\End_V(m,n))^{-*}
\subset W
\]
satisfy~(\ref{V}).%
\end{proof}

Let $\kappa \in C^*_\sfP(V;V) = \susp \Der(\sfM,\End_V)^{-*}$ and denote by
$u \in \Lin_{\Sigma-\Sigma}(E,\End_V)$ the restriction of $\desusp
\kappa$ to the space of
generators of $\sfM = \freePROP(E)$, $u:= \desusp \kappa |_E$. Let $U :
\sfM \to \End_V$ be the extension of $u$ into a \PROP{} homomorphism.

\begin{proposition}
\label{podzim}
Let $L_{\mbox{\o}} = (C^*_\sfP(V;V),h_2,h_3,\ldots)$ be the \li-structure
corresponding to the trivial $\sfP$-algebra. Under the above
notation,
one has the following equality of elements of
$\Lin_{\mbox {\scriptsize $\Sigma$-$\Sigma$}}(E,\End_V)$:
\begin{equation}
\label{je_to_trochu_lepsi}
U \circ \pa|_E = 
\raisebox{1.2em}{
\setlength{\unitlength}{1.8em}
\vector(0,-1){1}
}
\left.\left(
-\frac 1{2!} h_2(\kappa,\kappa) + \frac 1{3!}
h_3(\kappa,\kappa,\kappa) + \frac 1{4!}
h_4(\kappa,\kappa,\kappa,\kappa) - \cdots\right)\right|_E
\end{equation}
\end{proposition}

\begin{proof}
Since $L_{\mbox{\o}}$ corresponds to the trivial $\sfP$-structure, the map
$\beta$ in~(\ref{stanu_se_UL_instruktorem?}) is zero and the decorated
graph $\Gamma_{\{\beta\}}^{\{\Rada v1k\}}[u,\ldots,u]$ may be nontrivial
only if  $\Gamma$ has precisely $k$ vertices, all decorated by $u$. 
Therefore~(\ref{nabla}), with $f_1 = \cdots = f_k = \kappa$, describes
the $k$-th homogeneous component of the extension of $u$ into a 
homomorphism $U : \sfM \to \End_V$ composed with $\partial$. 
The signs in~(\ref{je_to_trochu_lepsi}) are induced by
(de)suspensions.%
\end{proof}

\section{Strongly homotopy algebras}
\label{sh}

In this short section we indicate how the methods of this
paper generalize to cohomology of strongly homotopy algebras and how
``curved'' $L_\infty$-algebras naturally arise in this context.

We will consider $L_\infty$-algebras $L=(W,l_0,l_1,l_2,\ldots)$ with
possibly nontrivial $l_0 \in W^1$.  These generalized
$L_\infty$-algebras can be defined by allowing $l_k$ for $k=0$ 
in~\cite[Definition~2.1]{lada-markl:CommAlg95}; we leave the
details for the reader. Axiom~(2) of~\cite{lada-markl:CommAlg95} for
$n=0$ gives $l_1 \circ l_0 = 0$ and for $n=1$
\begin{equation}
\label{pisu_ve_Dvore_Kralove}
0 = l_1(l_1 (w)) + l_2(l_0,w),\ w \in W.
\end{equation}
Therefore $l_1$ need not be a differential if $l_0 \not= 0$.  We will
call such an \li-algebra {\em curved\/} and $l_0$ its {\em
curvature\/}.  If $l_0 = 0$ we say that $L$ is {\em flat\/}. Flat
\li-algebras are thus ordinary \li-algebras without the $l_0$
term. \li-algebras with $l_0 = l_1 =0$ are sometimes called {\em
minimal\/}. The following statement is~\cite[Theorem~2.6.1]{merkulov},
slightly generalized by allowing curved $L_\infty$-algebras.

\begin{proposition}
\label{zavody}
Let $L= (W,l_0,l_1,l_2,\ldots)$ be an \li-algebra satisfying~(\ref{V})
and $\kappa \in W^1$ an \underline{arbitrar}y element. Then $L_\kappa
:=(W,l^\kappa_0,l^\kappa_1,l^\kappa_2,\ldots)$ with
\begin{eqnarray*}
l^\kappa_n(\Rada w1n) &:=& \sum_{s \geq 0} (-1)^{{sn} + {s+1\choose2}}
\frac 1{s!} 
l_{n+s}(\underbrace {\kappa,\ldots,\kappa}_{\mbox{\scriptsize 
$s$\ times}},\Rada w1n) 
\\
&&  \hskip -6em =
l_n(\Rada w1n) -(-1)^{n} 
l_{n+1}(\kappa,\Rada w1n) -\frac 12 l_{n+2}(\kappa,\kappa,\Rada w1n)
+ \cdots
\end{eqnarray*}
is an \li-algebra satisfying~(\ref{V}) whose curvature $l^\kappa_0$ equals
\[
l_0^\kappa = \sum_{s\geq 0} (-1)^{s+1 \choose 2} l_s(\rada \kappa\kappa) =
l_0 - l_1(\kappa) -\frac 1{2!}  l_2(\kappa,\kappa)
+ \frac 1{3!}  l_3(\kappa,\kappa,\kappa) +
 \cdots.
\]
\end{proposition}

The proof is a direct verification, see~\cite{merkulov}.
Let us remark that there is another sign convention for \li-algebras
used for example in~\cite{getzler:AM:09} related to the one introduced
in~\cite{lada-markl:CommAlg95} and used in this paper by $l_n
\leftrightarrow (-1)^{n+1\choose2}l_n$, $n \geq 0$. In this
convention, all terms in the above sums have the $+$ sign.

We will call $L_\kappa$ the {\em $\kappa$-twisting\/} of $L$.  Observe
that, if $L$ is flat and $\kappa$ satisfies the master
equation~(\ref{master1}) in $L$, then $L_\kappa$ is an ordinary flat
$L_\infty$-algebra. Proposition~\ref{zavody} then defines the
classical twisting of an $L_\infty$-algebra by a Maurer-Cartan
element, see for example~\cite[Lemma~4.4]{laan-defo}
or~\cite[Proposition~4.4]{getzler:AM:09}.  We will also need the
following elementary lemma whose proof is straightforward.

\begin{lemma}
\label{kousek_od_Zirce}
Suppose that, under assumptions of Proposition~\ref{zavody}, $W$ is
equipped with a degree $+1$ differential $d$ such that all operations
$l_n: \otexp Wn \to W$ 
of the algebra $L$ are chain maps (this, in particular, means
that $d l_0 = 0$). Suppose moreover that $\kappa \in W^1$ satisfies
\[
d(\kappa) = l_0 - l_1(\kappa) -\frac 1{2!}  l_2(\kappa,\kappa)
+ \frac 1{3!}  l_3(\kappa,\kappa,\kappa) +
\frac 1{4!}  k_4(\kappa,\kappa,\kappa,\kappa) - \cdots.
\]
Then $\overline{L}_\kappa:=(W,0,l^\kappa_1+d,l^\kappa_2,l^\kappa_3,\ldots)$, 
where $l^\kappa_n$ are as in
Proposition~\ref{zavody}, is a flat \li-algebra.
\end{lemma}

Let $\sfP$ be a $\bfk$-linear \PROP. 
Strongly homotopy $\sfP$-(bi)algebras are,
by definition, algebras
over the minimal model $(\sfM,\pa)$ of $\sfP$. 
Let $V_{\mbox{\o}}$ be the \sfP-(bi)algebra
whose all structure operations are trivial and 
$L_{\mbox{\o}} = (\CP^*(V,V),h_0=0,h_1 = 0,h_2,h_3,\ldots)$ the flat
\li-algebra constructed in Theorem~\ref{main} corresponding to
$V_{\mbox{\o}}$. The minimality $h_1 =0$
of $L_{\mbox{\o}}$ follows from the minimality of $(\sfM,\pa)$. 
Assume that $V$ is graded, with a degree $+1$ differential
$d$. Slightly abusing the notation, we will denote by $d$ also the
induced differential on $\End_V$. It is immediately clear that all
operations $l_n$ are chain maps. The following proposition is a
version of Lemma~5.11 of~\cite{laan-defo}.

\begin{proposition}
\label{pisu_v_aute_ve_Dvore}
There is a one-to-one correspondence between elements $\kappa \in
\CP^1(V,V)$ satisfying
\begin{equation}
\label{44}
d \kappa = -\frac 1{2!} h_2(\kappa,\kappa) + \frac 1{3!}
h_3(\kappa,\kappa,\kappa) + \frac 1{4!}
h_4(\kappa,\kappa,\kappa,\kappa) - \cdots
\end{equation}
in $L_{\mbox{\o}}$ and strongly homotopy \sfP-(bi)algebra structures on $V$.
\end{proposition}

\begin{proof}
Let $U :  \sfM \to \End_V$ be, as in
Proposition~\ref{podzim}, the (unique) homomorphism extending $\desusp
\kappa |_E$. Using~(\ref{je_to_trochu_lepsi}), one easily
sees that~(\ref{44}) is equivalent to $dU = U \partial$ which means
that $U$ is a homomorphism of dg-\PROP{s} defining a strongly homotopy
$\sfP$-algebra.%
\end{proof}

The interpretation of homotopy structures in terms of solutions of
the (generalized) Maurer-Cartan equation was found by van~der~Laan for
homotopy $\calP$-algebras~\cite{laan-defo}. The generalization to homotopy
(bi)algebras over a \PROP\ given here was independently found by
Merkulov-Vallette~\cite{merkulov-vallette}. 

We are finally ready to analyze the structure of the deformation
complex of a strongly homotopy (bi)algebra.
Let $\kappa \in \CP^1(V,V)$ be as in
Proposition~\ref{pisu_v_aute_ve_Dvore}. The curvature of the
$\kappa$-twisting 
\[
L_\kappa =
(\CP^*(V,V),h^\kappa_0,\delta_\kappa,h^\kappa_2,h^\kappa_3,\ldots)
\]
of $L_{\mbox{\o}}$ equals $d\kappa$. Assumptions of
Lemma~\ref{kousek_od_Zirce} are satisfied, therefore 
\[
\overline{L}_\kappa =(\CP^*(V,V),l_0 =0,d+
\delta_\kappa,h^\kappa_2,h^\kappa_3,\ldots)
\] 
is a flat \li-structure 
which induces a Lie bracket on the
cohomology 
\[
H_\sfP^*(B,B) := H^*(\CP^*(V,V), d+ \delta_\kappa)
\] 
of the
strongly homotopy (bi)algebra $B$ corresponding to $\kappa$.

\begin{example}
{\rm 
Let us illustrate the above analysis on \ai- (strongly homotopy
associative) algebras~\cite{stasheff:TAMS63}.  They are algebras over
the minimal model $\Ass_\infty$ of the operad $\Ass$ for associative
algebras. It immediately follows from the description of $\Ass_\infty$
given for instance in~\cite[Example~4.8]{markl:zebrulka} that
\[
C^*_{\ssAss}(V,V) = \prod_{n \geq 2}\Lin^{n-*-1}(\otexp Vn,V).
\]
The only nontrivial operation of the flat algebra $L_{\mbox{\o}}$ is the bilinear
bracket $h_2$ given by
\begin{equation}
\label{znechuceni}
h_2\left((\phi_2,\phi_3,\ldots),(\psi_2,\psi_3,\ldots)\right)_n 
:=  \sum_{i+j = n+1} [\phi_i,\psi_j],
\end{equation}
where $\phi_s,\psi_s \in \Lin(\otexp Vs,V)$, $s \geq 2$,
the subscript $n$ denotes the component in $\Lin(\otexp Vn,V)$
and $[-,-]$ is the Gerstenhaber bracket of Hochschild
cochains~\cite{gerstenhaber:AM63}. 

Equation~(\ref{44}) for $\kappa = (\mu_2,\mu_3,\ldots) \in
C^1_{\ssAss}(V,V) = \prod_{n \geq 2}\Lin^{n-2}(\otexp Vn,V)$, 
expanded into homogeneous components,
reads
\[
0 = d \mu_n + \frac 12 \sum_{i+j = n+1} [\mu_i,\mu_j],\ n \geq 2,
\]
which we easily recognize as the axiom for \ai-algebras in the
form~\cite[Section~1.4]{markl:JPAA92}.

The $\kappa$-twisting $L_\kappa$ of $L_{\mbox{\o}}$ equals $L_\kappa =
(C^*_{\ssAss}(V,V),l_0,\delta_\kappa,l_2,0,0,\ldots)$, with the curvature
\[
l_0 = (d\mu_2,d\mu_3,d\mu_4,\ldots),
\] 
$\delta_\kappa$ given by
\[
\delta(f_2,f_3,\ldots)_n := \sum_{i+j = n+1} [\mu_i,f_j]
\]
and $l_2 = h_2$ as in~(\ref{znechuceni}).  In the ``flattened'' algebra
\[
\overline{L}_\kappa = (C^*_{\ssAss}(V,V),l_0 = 0, d+
\delta_\kappa,l_2,0,0,\ldots),
\] 
$d + \delta_\kappa$ is the
differential on the cochain complex defining the cohomology of the
\ai-algebra determined by $\kappa$ with coefficients in
itself~\cite[Section~2.2]{markl:JPAA92}.  
}\end{example}

\section{The Gerstenhaber-Schack cohomology of bialgebras}
\label{GS}\label{3}

In this section we show how to construct, by applying methods of
Section~\ref{1bis} to the differential of the minimal model
$(\sfM_\sfB,\psfb)$ for the bialgebra \PROP\ $\sfB$, an explicit Lie
bracket on the Gerstenhaber-Schack cohomology of a bialgebra.
Formulas for the differential $\psfb$ of $\sfM_\sfB$ are given 
in~\cite[Eqn.~3.1]{umble-saneblidze:KK}.

Recall that a ($\Ass$-)bialgebra $B$ is a vector space $V$ with
a {\em multiplication\/} $\mu : V \ot V \to V$ and a {\em
comultiplication\/} (also called a {\em diagonal\/}) $\Delta : V \to
V\ot V$. The multiplication is associative:
\[
\mu(\mu \ot \iden_V) = \mu(\iden_V \ot \mu),
\]
the comultiplication is coassociative:
\[
(\iden_V \ot \Delta)\Delta = (\Delta \ot \iden_V)\Delta
\]
and the usual compatibility relation between $\mu$ and $\Delta$ is
assumed:
\begin{equation}
\label{compatibility}
\Delta \circ \mu = (\mu \ot \mu) T_{\sigma(2,2)} (\Delta \ot \Delta),
\end{equation}
where $T_{\sigma(2,2)}: V^{\ot 4} \to V^{\ot 4}$ is defined by
\[
T_{\sigma(2,2)}(v_1 \ot v_2 \ot v_3 \ot v_4) := v_1 \ot v_3 \ot v_2
\ot v_4,
\]
for $v_1,v_2,v_3,v_4 \in V$. Compatibility~(\ref{compatibility})
of course expresses
the fact that
\[
\Delta(u \cdot v) = \Delta(u) \cdot \Delta(v),\ u,v \in V,
\]
where $u \cdot v :=\mu(u,v)$ and the dot $\ \cdot\ $ in the right hand
side denotes the multiplication induced on $V \ot V$ by $\mu$.

Let us recall that in the definition of the Gerstenhaber-Schack
cohomology~\cite{gerstenhaber-schack:Proc.Nat.Acad.Sci.USA90} one
considers the bigraded vector space
\[
C^{*,*}_\GS (B;B) := \bigoplus_{p,q \geq 1} C^{p,q}_\GS (B;B),
\]
where
\[
C^{p,q}_\GS (B;B) := \Lin(\otexp Vp,\otexp Vq).
\]
It will be useful to introduce the {\em biarity\/} of a function $f
\in C^{p,q}_\GS (B;B)$ as the couple $\biar(f) := (q,p)$.  For each $q
\geq 2$, the iterated diagonal
\[
\Delta^{[q]}:= (\Delta \ot
{\iden}^{\ot (q-2)})\circ(\Delta \ot
{\iden}^{\ot (q-1)}) \circ \cdots \circ \Delta :  V \to \otexp Vq
\]
induces on $\otexp Vq$ the structure of a $(V,\mu)$-bimodule, by
\begin{eqnarray*}
u(\otRada v1q)&: = &\Delta^{[q]}(u) \cdot (\otRada v1q),\ \mbox { and }
\\
(\otRada v1q)u &: = & (\otRada v1q) \cdot \Delta^{[q]}(u),
\end{eqnarray*}
where $\ \cdot\ $ denotes the multiplication induced on $\otexp Vq$ by
$\mu$.  Therefore it makes sense to define
\[
d_1 : C^{p,q}_\GS (B;B) \to C^{p+1,q}_\GS (B;B)
\]
to be the Hochschild differential of the algebra $(V,\mu)$ with
coefficients in the $(V,\mu)$-bimodule $\otexp Vq$. The ``coHochschild''
differential
\[
d_2 : C^{p,q}_\GS (B;B) \to C^{p,q+1}_\GS (B;B)
\]
is defined in dual manner. It turns out that
$(C^{*,*}_\GS (B;B), d_1 + d_2)$
forms a bicomplex shown in Figure~\ref{fig1}.
\begin{figure}[t]
\setlength{\unitlength}{.8em}
\begin{picture}(20,15)(-3,0)
\put(0,1){
\put(2.5,0){
\put(0,0){\makebox(0,0){$\CGSbi 11$}}
\put(10,0){\makebox(0,0){$\CGSbi 21$}}
\put(20,0){\makebox(0,0){$\CGSbi 31$}}
\put(30,0){\makebox(0,0){$\CGSbi 41$}}
}
\put(2.5,5){
\put(0,0){\makebox(0,0){$\CGSbi 12$}}
\put(10,0){\makebox(0,0){$\CGSbi 22$}}
\put(20,0){\makebox(0,0){$\CGSbi 32$}}
\put(30,0){\makebox(0,0){$\CGSbi 42$}}
}
\put(2.5,10){
\put(0,0){\makebox(0,0){$\CGSbi 13$}}
\put(10,0){\makebox(0,0){$\CGSbi 23$}}
\put(20,0){\makebox(0,0){$\CGSbi 33$}}
\put(30,0){\makebox(0,0){$\CGSbi 43$}}
}
\put(5.5,0){\multiput(0,0)(10,0){4}{
            \vector(1,0){3}\put(-1.2,1){\makebox(0,0){$d_1$}}}}
\put(5.5,5){\multiput(0,0)(10,0){4}{
            \vector(1,0){3}\put(-1.2,1){\makebox(0,0){$d_1$}}}}
\put(5.5,10){\multiput(0,0)(10,0){4}{
            \vector(1,0){3}\put(-1.2,1){\makebox(0,0){$d_1$}}}}
\put(2.5,1){\multiput(0,0)(0,5){3}{
            \vector(0,1){3}\put(-1,1.2){\makebox(0,0){$d_2$}}}}
\put(12.5,1){\multiput(0,0)(0,5){3}{
            \vector(0,1){3}\put(-1,1.2){\makebox(0,0){$d_2$}}}}
\put(22.5,1){\multiput(0,0)(0,5){3}{
            \vector(0,1){3}\put(-1,1.2){\makebox(0,0){$d_2$}}}}
\put(32.5,1){\multiput(0,0)(0,5){3}{
            \vector(0,1){3}\put(-1,1.2){\makebox(0,0){$d_2$}}}}
}
\end{picture}
\caption{\label{fig1}The Gerstenhaber-Schack bicomplex.}
\end{figure}
The {\em Gerstenhaber-Schack\/} cohomology of $B$ with coefficients in
$B$ is the cohomology of its (regraded) total complex
\[
H^*_\GS(B;B) := H^*(C^*_\GS(B;B),d_\GS),
\]
where
\[
C^*_\GS(B;B) := \bigoplus_{* = p+q-2} C^{p,q}_\GS (B;B)\
\mbox { and }\
d_\GS := d_1 + d_2.
\]

Let us compare now the Gerstenhaber-Schack cohomology with the cohomology
$H^*_\sfB(B;B)$ recalled in~(\ref{Tel-Aviv}), 
where $\sfB$ is the \PROP\ for bialgebras. To this
end, we need to review some facts about the minimal model of $\sfB$.
For a generator $\xi^m_n$ of biarity $(m,n)$ ($n$ `inputs' and $m$
`outputs,' $m,n \geq 0$), let $\sigmaspan(\xi) : = \bfk[\Sigma_m]
\otimes \bfk \cdot \xi \otimes \bfk[\Sigma_n]$, with the obvious
mutually compatible left $\Sigma_m$- right $\Sigma_n$-actions.
The minimal model $\sfM_\sfB = \sfM$ of
$\sfB$ is of the form $\sfM = (\freePROP(\Xi),\psfb)$, where $\Xi :=
\sigmaspan(\{\xi^m_n\}_{m,n \in I})$ with
\begin{equation}
\label{SM_je_hovado}
I := \{ m,n \geq 1,\ (m,n) \not= (1,1)\},
\end{equation}
see~\cite[Theorem~16]{markl:ba}.  The differential $\psfb$ of
$\sfM_\sfB$ is explicitly determined by~\cite[Eqn.~3.1]{umble-saneblidze:KK}.

While general free \PROP{s} are complicated objects, the fact that
the $\Sigma$-bimodule $\Xi$ is a~direct sum of regular
representations implies the following relatively simple description of
$\freePROP(\Xi)$.
Let us agree that, in this section, 
by directed a graph we mean a finite, not necessary
connected, graph $G$ such that
\begin{itemize}
\item[(i)]
each edge $e \in \edge(G)$ 
is equipped with a direction and there are no directed cycles (wheels)
in $G$.
\end{itemize}
The direction of edges determines at
each vertex $v \in \vert(G)$ of a directed graph $G$ a disjoint decomposition
\[
\edge(v) = \In(v) \sqcup \Out(v)
\]
of the set of edges adjacent to $v$ into the set $\In(v)$ of incoming
edges and the set $\Out(v)$ of outgoing edges.  
We will also assume that
\begin{itemize}
\item[(ii)]
for each vertex $v \in \vert(G)$, the sets $\In(v)$ and $\Out(v)$ are
linearly ordered. 
\end{itemize}

The pair $\biar(v):= (\#(\Out(v)),\#(\In(v)))$ is called the {\em
biarity\/} of $v$. A vertex $v \in \Vert(G)$ is {\em
binary\/} if $\biar(v) = (1,2) \mbox { or } (2,1)$.  
Next assumption we impose on the graphs in this section is that
\begin{itemize}
\item[(iii)] 
$G$ has no vertices of biarity $(0,0)$ or $(1,1)$.
\end{itemize}
Vertices of biarity $(1,0)$ are called the {\em input vertices\/} and
vertices of biarity $(0,1)$ the {\em output vertices\/} of $G$.
We finally denote by $\ggg mn$ the set of isomorphism classes of 
directed graphs $G$ satisfying (i)--(iii) above such that
\begin{itemize}
\item[(iv)] 
the input vertices of $G$ are labeled by
$\{1, \dots, n\}$ and the output vertices by $\{1, \dots, m\}$.
\end{itemize}
With this notation,
\[
\freePROP(\Xi)(m,n) \cong \Span(\{G\}_{G \in \ggg mn}),\ m,n \geq 0.
\]

The degree $-1$ differential $\psfb$ on $\freePROP(\Xi)$ is of course
uniquely determined by its values $\psfb (\xi^m_n) \in \freePROP(\Xi)(m,n)$,
$(m,n) \in I$, on the generators of $\Xi$. We will see more explicitly
below how these values look, now just write
\begin{equation}
\label{je_to_takova_glegla}
\psfb(\xi^m_n) = \sum_{s \in S^m_n} \epsilon^s \cdot G_s,
\end{equation}
where  $G_s \in \ggg mn$, $S^m_n$ is a finite indexing set depending
on $(m,n) \in I$, and $\epsilon \in \bfk$. 

It follows from the results of~\cite[Section~4]{umble-saneblidze:KK}
that the minimal model of $\sfB$ is the cellular chain complex of a
sequence of finite polytopes whose faces are indexed  
by directed graphs that we may
in~(\ref{je_to_takova_glegla}) assume $\epsilon^s \in
\{-1,1\}$. In particular, the minimal model of $\sfB$ is defined over
the integers! This has to be compared to a similar argument proving
the integrality of the minimal model of the operad for associative
algebras based on the existence of the
associahedra~\cite[Example~4.8]{markl:zebrulka}. 

In the same manner, any derivation $F \in
\Der(\sfM,\End_V)$ is uniquely determined by specifying, for each $(m,n)
\in I$, multilinear maps
\[
F(\xi^m_n) \in \End_V(m,n) = \Lin(\otexp
Vn,\otexp Vm) = \CGSbi nm .
\]
This defines an isomorphism
\[
C_\sfB^*(V;V) = \hskip .3 em\susp \Der(\sfM,\End_V)^{-*} \cong  C^*_\GS(B;B)
\]
which we use to identify $C_\sfB^*(V;V)$ with $ C^*_\GS(B;B)$.  Let us
inspect how the differential $\delta_\sfB$ of the cochain complex
$C_\sfB^*(V;V)$ acts on some $f \in \CGSbi pq$. If we denote by
$(\delta_\sfB f)^m_n$ the component of $\delta_\sfB f$ in $\CGSbi nm$, then
\begin{equation}
\label{psano_v_Koline}
(\delta_\sfB f)^m_n = \sum_{s \in S^m_n}
\sum_{v \in \Vert(G_s)} \epsilon^s \cdot G^v_s[f],
\end{equation}
where $G^v_s[f]$ is an element of $\Lin(\otexp
Vn,\otexp Vm)$ which is nontrivial only if the vertex $v \in
\Vert(G_s)$ has biarity $(q,p)$ and if all remaining vertices of
$G_s$ are binary. In this case $G^v_s[f]$ is obtained by
decorating the vertex $v$ with $f$, vertices of biarity $(1,2)$
with the multiplication $\mu$, vertices of biarity $(2,1)$ with
the comultiplication $\Delta$, and then performing the
compositions indicated by the graph $G_s$.  We leave as an
exercise to show that formula~(\ref{psano_v_Koline}) indeed
follows from the definition of $\delta_\sfP$ with $\sfP = \sfB$
recalled in Section~\ref{intro}.

It turns out (see Appendix~B on page~\pageref{apB}) 
that $(\delta_\sfB f)^m_n \not= 0$ only if $(m,n) = (q,p+1)$
or $(m,n) = (q+1,p)$ and that
\begin{equation}
\label{Merkulov_je_rusky_hovado}
(\delta_\sfB f)^q_{p+1} = d_1 f\ \mbox { and }\
(\delta_\sfB f)^{q+1}_p = d_2 f.
\end{equation}
This shows that the cochain complexes $(C^*_\GS(B;B),d_\GS)$ and
$(C^*_\sfB(B;B),\delta_\sfB)$ are isomorphic and that our cohomology
$H^*_\sfB(B;B)$  coincides with the Gerstenhaber-Schack
cohomology of $B$.

The $L_\infty$-structure on $C^*_\GS(B;B)$ announced in the Abstract
is determined by an obvious generalization
of~(\ref{psano_v_Koline}).  The components $l_k(\Rada f1k)^m_n \in
\CGSbi mn$ of the bracket $l_k(\Rada f1k)$ can be computed as
\begin{eqnarray}
\label{dnes_prijede_ten_American}
\lefteqn{
l_k(\Rada f1k)^m_n = \hskip 5em}
\\
&& \hskip 5em \nonumber  
(-1)^{\nu(\Rada f1n)} \cdot
 \sum_{s \in S^m_n}\sum_{\Rada v1k \in \Vert(G_s)}
\epsilon^s G_s^{\{\Rada v1k\}}[\Rada f1k],
\end{eqnarray}
where $\nu(\Rada f1n)$ is as in~(\ref{dnes_jsem_byl_na_sipcich})
and in the second summation we assume that all vertices $\Rada
v1k$ are mutually different.  The multilinear map $G_s^{\{\Rada
v1k\}}[\Rada f1k] \in \Lin(\otexp Vn,\otexp Vm)$
is nonzero only if $\biar(f_i) = \biar(v_i)$ for $1 \leq i \leq
k$ and if all other remaining vertices of $G_s$ are binary. When
it is so, then $G_s^{\{\Rada v1k\}}[\Rada f1k]$ is defined by
decorating the vertices $v_i$ with $f_i$, $1 \leq i \leq k$,
vertices of biarity $(1,2)$ with $\mu$, vertices of biarity
$(2,1)$ with $\Delta$, and then performing the compositions
indicated by the graph $G_s$.

Before we give examples of these $L_\infty$-brackets, we
need to recall the calculus of fractions introduced
in~\cite{markl:ba}. A fraction is a special type of a composition of
elements of a \PROP\ defined using a restricted class of permutations
which we we need to recall first.

For $k ,l \geq 1$ and $1 \leq i \leq kl$, let $\sigma(k,l) \in
\Sigma_{kl}$ be the permutation given~by
\[
\sigma(i) := k(i-1 - (s-1)l) + s,
\]
where $s$ is such that $(s-1)l < i \leq sl$. Permutations of
this form are called {\em special permutations\/}.
An example is the permutation $\sigma(2,2)$
in~(\ref{compatibility}). Another example is
\[
\sigma(3,2) :=
\left(
\begin{array}{cccccc}
1 & 2 & 3 & 4 & 5 & 6
\\
1 & 4& 2 & 5 & 3 & 6
\end{array}
\right)
= \hskip 3mm
{
\unitlength=1.000000pt
\begin{picture}(50.00,25.00)(0.00,12.00)
\put(50.00,0.00){\makebox(0.00,0.00){\scriptsize$\bullet$}}
\put(50.00,30.00){\makebox(0.00,0.00){\scriptsize$\bullet$}}
\put(40.00,0.00){\makebox(0.00,0.00){\scriptsize$\bullet$}}
\put(40.00,30.00){\makebox(0.00,0.00){\scriptsize$\bullet$}}
\put(30.00,0.00){\makebox(0.00,0.00){\scriptsize$\bullet$}}
\put(30.00,30.00){\makebox(0.00,0.00){\scriptsize$\bullet$}}
\put(20.00,0.00){\makebox(0.00,0.00){\scriptsize$\bullet$}}
\put(20.00,30.00){\makebox(0.00,0.00){\scriptsize$\bullet$}}
\put(10.00,0.00){\makebox(0.00,0.00){\scriptsize$\bullet$}}
\put(10.00,30.00){\makebox(0.00,0.00){\scriptsize$\bullet$}}
\put(0.00,0.00){\makebox(0.00,0.00){\scriptsize$\bullet$}}
\put(0.00,30.00){\makebox(0.00,0.00){\scriptsize$\bullet$}}
\put(0.00,30.00){\line(0,1){0.00}}
\put(50.00,20.00){\line(0,-1){10.00}}
\put(40.00,20.00){\line(-1,-1){10.00}}
\put(30.00,20.00){\line(-2,-1){20.00}}
\put(20.00,20.00){\line(2,-1){20.00}}
\put(10.00,20.00){\line(1,-1){10.00}}
\put(0.00,20.00){\line(0,-1){10.00}}
\end{picture}} \hskip 2mm
\]

\noindent
with the convention that the `flow diagrams'
should be read from the bottom to the top.

Let $\sfP$ be a \PROP.
Let $k,l \geq 1$, $\Rada a1l \geq 1$, $\Rada b1k \geq 1$,
$\Rada A1l \in \sfP(a_i,k)$ and $\Rada B1k \in \sfP(l,b_j)$.
Then the {\em $(k,l)$-fraction\/} is defined as
\[
\frac{A_1 \cdots A_l}{B_1 \cdots B_k} :=
(A_1 \otimes \cdots \otimes A_l)
\circ \sigma(k,l) \circ
(B_1 \otimes \cdots \otimes B_k)
\in \sfP(a_1+\cdots + a_l,b_1+\cdots + b_k).
\]

If $k=1$ or $l=1$, the $(k,l)$-fractions are just the `operadic'
compositions:
\[
\frac{A_1 \ot \cdots \ot A_l}{B_1} = (A_1 \ot \cdots \ot A_l) \circ B_1
\ \mbox { and }\
\frac{A_1}{B_1 \ot \cdots \ot B_k} = A_1 \circ (B_1 \ot \cdots \ot B_k).
\]
As a more complicated example,
consider $\zeroatwo a, \zeroatwo b \in \sfP(*,2)$ and $\twoazero c,
\twoazero d \in \sfP(2,*)$. Then

\vglue -1.3em
\[
\frac{\zeroatwo a \hskip .2em  \zeroatwo b}%
     {\twoazero c \hskip .2em \twoazero d}
=
(\raisebox{-2pt}{\zeroatwo a} \hskip -.3em
\ot \raisebox{-2pt}{\zeroatwo b}
\hskip -4pt)
\circ \sigma(2,2) \circ
(\raisebox{-2pt}{\twoazero c}
\hskip -.3em\ot \raisebox{-2pt}{\twoazero d} \hskip -4pt)
=
{
\unitlength=.5pt
\begin{picture}(70.00,70.00)(0.00,30.00)
\put(50.00,10.00){\makebox(0.00,0.00)[l]{$d$}}
\put(10.00,10.00){\makebox(0.00,0.00)[l]{$c$}}
\put(50.00,60.00){\makebox(0.00,0.00)[l]{$b$}}
\put(10.00,60.00){\makebox(0.00,0.00)[l]{$a$}}
\put(60.00,50.00){\line(0,-1){30.00}}
\put(40.00,0.00){\line(0,1){20.00}}
\put(70.00,0.00){\line(-1,0){30.00}}
\put(70.00,20.00){\line(0,-1){20.00}}
\put(40.00,20.00){\line(1,0){30.00}}
\put(40.00,70.00){\line(0,-1){20.00}}
\put(70.00,70.00){\line(-1,0){30.00}}
\put(70.00,50.00){\line(0,1){20.00}}
\put(40.00,50.00){\line(1,0){30.00}}
\put(50.00,50.00){\line(-1,-1){30.00}}
\put(20.00,50.00){\line(1,-1){30.00}}
\put(0.00,0.00){\line(0,1){20.00}}
\put(30.00,0.00){\line(-1,0){30.00}}
\put(30.00,20.00){\line(0,-1){20.00}}
\put(0.00,20.00){\line(1,0){30.00}}
\put(10.00,40.00){\line(0,-1){20.00}}
\put(10.00,50.00){\line(0,-1){10.00}}
\put(30.00,70.00){\line(-1,0){30.00}}
\put(30.00,50.00){\line(0,1){20.00}}
\put(0.00,50.00){\line(1,0){30.00}}
\put(0.00,70.00){\line(0,-1){20.00}}
\end{picture}}\hskip 2mm.
\]
Similarly, for $\zeroathree x, \zeroathree y \in \sfP(*,3)$ and
$\twoazero z, \twoazero u, \twoazero v \in \sfP(2,*)$,

\vglue -1.3em
\[
\frac{\zeroathree x \hskip 4mm \zeroathree y}%
     {\twoazero z \hskip 1mm \twoazero u \hskip 1mm \twoazero v}
=
(\hskip -2pt \raisebox{-2pt}{\zeroathree x}\hskip -.2em \ot\hskip -.2em
\raisebox{-2pt}{\zeroathree y}\hskip -2pt)
\circ \sigma(3,2) \circ
(\raisebox{-2pt}{\twoazero z}\hskip -.2em
\ot  \raisebox{-2pt}{\twoazero u}  \hskip -.2em
\ot \raisebox{-2pt}{\twoazero v}\hskip -2pt)
=
{
\unitlength=.5pt
\begin{picture}(110.00,70.00)(0.00,30.00)
\put(90.00,10.00){\makebox(0.00,0.00)[l]{$v$}}
\put(50.00,10.00){\makebox(0.00,0.00)[l]{$u$}}
\put(10.00,10.00){\makebox(0.00,0.00)[l]{$z$}}
\put(90.00,60.00){\makebox(0.00,0.00){$y$}}
\put(20.00,60.00){\makebox(0.00,0.00){$x$}}
\put(110.00,20.00){\line(-1,0){30.00}}
\put(110.00,0.00){\line(0,1){20.00}}
\put(80.00,0.00){\line(1,0){30.00}}
\put(80.00,20.00){\line(0,-1){20.00}}
\put(40.00,0.00){\line(0,1){20.00}}
\put(70.00,0.00){\line(-1,0){30.00}}
\put(70.00,20.00){\line(0,-1){20.00}}
\put(40.00,20.00){\line(1,0){30.00}}
\put(0.00,0.00){\line(0,1){20.00}}
\put(30.00,0.00){\line(-1,0){30.00}}
\put(30.00,20.00){\line(0,-1){20.00}}
\put(0.00,20.00){\line(1,0){30.00}}
\put(60.00,20.00){\line(1,1){30.00}}
\put(20.00,20.00){\line(2,1){60.00}}
\put(30.00,50.00){\line(2,-1){60.00}}
\put(20.00,50.00){\line(1,-1){30.00}}
\put(100.00,50.00){\line(0,-1){30.00}}
\put(70.00,70.00){\line(0,-1){20.00}}
\put(110.00,70.00){\line(-1,0){40.00}}
\put(110.00,50.00){\line(0,1){20.00}}
\put(70.00,50.00){\line(1,0){40.00}}
\put(0.00,70.00){\line(0,-1){20.00}}
\put(40.00,70.00){\line(-1,0){40.00}}
\put(40.00,50.00){\line(0,1){20.00}}
\put(0.00,50.00){\line(1,0){40.00}}
\put(10.00,50.00){\line(0,-1){30.00}}
\end{picture}}\hskip 2pt.
\]

To use formula~(\ref{dnes_prijede_ten_American}), we need of course to know the
differential $\psfb$ in the minimal model $\sfM =
(\freePROP(\Xi),\psfb)$ which determines the graphs $G_s$.
In~\cite{markl:ba}, $\psfb$ was described as a perturbation, $\psfb = \pa_0
+ \papert$, of its \hPROP{}-part $\pa_0$. Let us therefore give some
formulas for the unperturbed part $\pa_0$ first.

If we denote $\xi^1_2 = \jednadva$ and $\xi^2_1 = \dvajedna$, then
$\pa_0(\jednadva) = \pa_0(\dvajedna) =0$.  If $\xi^2_2 = \dvadva$,
then $ \pa_0(\dvadva) = \dvojiteypsilon$.  With the obvious, similar
notation,
\begin{eqnarray*}
\pa_0(\jednatri) &=& \ZbbZb - \bZbbZ,
\\
\label{jedna-ctyri}
\pa_0(\jednactyri) &=& \ZbbZbb - \bZbbZb + \bbZbbZ - \ZbbbZb - \bZbbbZ,
\\
\nonumber
\pa_0(\trijedna) &=& \ZvvZv - \vZvvZ,
\\
\nonumber
\pa_0(\dvatri) &=&
\dvacarkatri - \dvaZbbZb + \dvabZbbZ,
\\
\nonumber
\pa_0(\tridva) &=&  -
{
\unitlength=.2pt
\begin{picture}(40.00,50.00)(0.00,-50.00)
\put(20.00,-30.00){\line(0,1){10.00}}
\put(20.00,0.00){\line(0,-1){20}}
\bezier{20}(20.00,-20.00)(30.00,-10.00)(40.00,0.00)
\bezier{20}(20.00,-20.00)(10.00,-10.00)(0.00,0.00)
\bezier{20}(20.00,-30.00)(30.00,-40.00)(40.00,-50.00)
\bezier{20}(0.00,-50.00)(10.00,-40.00)(20.00,-30.00)
\end{picture}}
+ \ZvvZvdva - \vZvvZdva,\ \mbox { \&c.}
\end{eqnarray*}

And here is the differential $\psfb$ in its full beauty: $\psfb(\jednadva) =
0$, $\psfb(\dvajedna) = 0$, $\psfb(\jednatri) = \pa_0(\jednatri)$,
$\psfb(\trijedna) = \pa_0(\trijedna)$,
\def\mfrac#1#2{\frac{\raisebox{-.2em}{\hskip .2em#1}}%
{\raisebox{-.1em}{\hskip .2em#2}}}
\begin{eqnarray*}
\psfb(\dvadva) &=& \pa_0(\dvadva)
-
\mfrac{\jednadva \jednadva}{\dvajedna  \dvajedna},
\\
\psfb(\dvatri) &=& \pa_0(\dvatri)
+
\mfrac{\jednadva \jednadva}{\dvajedna  \dvadva}
-
\mfrac{\jednadva  \jednadva}{\dvadva  \dvajedna}
-
\mfrac{\bZbbZ \hskip 1em \jednatri}%
     {\dvajedna  \dvajedna  \dvajedna}
-
\mfrac{\hskip -.2em\jednatri \hskip 1em  \ZbbZb}{\dvajedna  \dvajedna  \dvajedna},
\\
\psfb(\gen 32) &=& \pa_0(\tridva)
-
\mfrac{\jednadva \dvadva}{\dvajedna \dvajedna}
+
\mfrac{\dvadva \jednadva}{\dvajedna  \dvajedna}
+
\mfrac{\jednadva \jednadva \jednadva}
     {\vZvvZ \hskip 1em \trijedna}
+
\mfrac{\jednadva \jednadva \jednadva}
     {\hskip -.2em\trijedna  \hskip 1em \ZvvZv},\ \mbox { \&c.}
\end{eqnarray*}

Let us finally see what formula~(\ref{dnes_prijede_ten_American}) tells us in
some concrete situations.  For $f,g \in \CGSbi 21$, the component
$l_2(f,g)^1_3$ is calculated as
\[
l_2(f,g)^1_3 =
(-1)^{|f|} \sum_{u,v} \psfb(\jednatri)^{\{u,v\}} [f,g]
=
(-1)^{|f|} \sum_{u,v}
\left(\ZbbZb - \bZbbZ\right)^{\{u,v\}} [f,g],
\]
where $|f| = 1$. Expanding the right-hand side, we obtain
\[
l_2(f,g)^1_3 = \raisebox{-.4em}{\bZbbZdec fg} +
\raisebox{-.4em}{\bZbbZdec gf} -
\raisebox{-.4em}{\ZbbZbdec fg} -
\raisebox{-.4em}{\ZbbZbdec gf},
\]
which we easily recognize as the ``classical'' Gerstenhaber bracket
$[f,g]$ of bilinear cochains.  The component $l_2(f,g)^2_2$ is given
by
\begin{eqnarray*}
l_2(f,g)^2_2 &=&
(-1)^{|f|} \sum_{u,v} \psfb(\dvadva)^{\{u,v\}} [f,g]
=
(-1)^{|f|} \sum_{u,v}
\left(\dvojiteypsilon  - \mfrac{\jednadva \hskip .2em
\jednadva}{\dvajedna \hskip .2em \dvajedna} \right)^{\{u,v\}} [f,g]
\\
&=&
 -\sum_{u,v}
\dvojiteypsilon^{\{u,v\}} [f,g] +
\sum_{u,v}
\mfrac{\jednadva \hskip .2em
\jednadva}{\dvajedna \hskip .2em \dvajedna}^{\{u,v\}} [f,g].
\end{eqnarray*}
The first term of the last line is zero, while the second one expands to
\[
l_2\left(\raisebox{-.4em}{$\jednadvadec f$},
   \raisebox{-.4em}{\jednadvadec g}   \right)^2_2 =
\frac{\jednadvadec f  \jednadvadec g}
     {\rule{0em}{1.8em}\dvajednadec \Delta \dvajednadec \Delta}
+
\frac{\jednadvadec g  \jednadvadec f}
     {\rule{0em}{1.8em}\dvajednadec \Delta \dvajednadec \Delta} \hskip
   .2em.
\]
This provides the first example of a ``non-classical''
bracket. In terms of elements,
\[
l_2(f,g)^2_2 (u,v) = f(u_{(1)} \ot v_{(1)}) \otimes g(u_{(2)} \ot v_{(2)}) +
g(u_{(1)} \ot v_{(1)})\otimes f(u_{(2)} \ot v_{(2)})
\]
for $\ u\otimes v \in \otexp V2$, where the standard Sweedler
notation for the coproduct is used.  An example of a triple bracket is
given by
\[
l_3\left(\raisebox{-.4em}{$\jednadvadec f$},
\raisebox{-.4em}{$\jednadvadec g$}, \raisebox{-.4em}{$\dvajednadec
h$}\right)^2_2 = \hskip .2em
\frac{\jednadvadec f  \jednadvadec g}
     {\rule{0em}{1.8em}\dvajednadec h \dvajednadec \Delta}
+ \frac{\jednadvadec g  \jednadvadec f}
     {\rule{0em}{1.8em}\dvajednadec h \dvajednadec \Delta}
+\frac{\jednadvadec f  \jednadvadec g}
     {\rule{0em}{1.8em}\dvajednadec \Delta \dvajednadec h}
+\frac{\jednadvadec g  \jednadvadec f}
     {\rule{0em}{1.8em}\dvajednadec \Delta \dvajednadec h} \hskip .2em.
\]
Slightly more complicated is:
\[
l_3\left(\raisebox{-.4em}{\jednadvadec f},
\raisebox{-.4em}{\jednadvadec g},
\raisebox{-.4em}{\jednatridec h}\right)^2_3 =
 \hskip .3em \frac{\bZbbZdec fg \hskip 1em \jednatridec h}
{\rule{0em}{1.8em}\dvajednadec \Delta \dvajednadec \Delta \dvajednadec
\Delta}
+\frac{\bZbbZdec gf \hskip 1em \jednatridec h}
{\rule{0em}{1.8em}\dvajednadec \Delta \dvajednadec \Delta \dvajednadec
\Delta}
+\frac{\jednatridec h  \hskip 1em \ZbbZbdec fg}
{\rule{0em}{1.8em}\dvajednadec \Delta \dvajednadec \Delta \dvajednadec
\Delta}
+\frac{\jednatridec h  \hskip 1em \ZbbZbdec gf}
{\rule{0em}{1.8em}\dvajednadec \Delta \dvajednadec \Delta \dvajednadec
\Delta}
\hskip .5em .
\]

We believe that the reader already understands the ideas behind our
definitions and that he or she can easily construct other 
examples of brackets using explicit formulas for the differential
$\psfb$~\cite[Eqn.~3.1]{umble-saneblidze:KK}.

Notice that the bottom row $(C^{*,1}_\GS(B,B),d_1)$ of the bicomplex
in Figure~\ref{fig1} is the Hochschild complex of the algebra $B =
(V,\mu)$ with coefficients in itself. For arbitrary $f_1,f_2 \in
C^{*,1}_\GS(B,B)$, the component $l_2(f_1,f_2)^1_*$ of $l_2(f_1,f_2)$
coincides with the classical Gerstenhaber bracket of the Hochschild
cochains~\cite{gerstenhaber:AM63} while the components $l_n(\Rada
f1n)^1_*$, $n \geq 3$, of the higher brackets are trivial. In this sense
our construction extends the classical Gerstenhaber bracket of
Hochschild cochains.

\begin{example}
{\rm
Let $V_{\mbox{\o}} := (V,\mu =0,\Delta = 0)$ be a 
vector space $V$ considered as
a bialgebra with trivial product and coproduct. Let $m \in
C^{2,1}_\GS(V_{\mbox{\o}},V_{\mbox{\o}}) = \Lin(\otexp V2,V)$,   $c \in
C^{1,2}_\GS(V_{\mbox{\o}},V_{\mbox{\o}}) = \Lin(V,\otexp V2)$ and $\kappa := m+c$. Finally,
let $(h_1,h_2,h_3,\ldots)$ be the $L_\infty$-structure on the
Gerstenhaber-Schack complex $C^*_\GS(V_{\mbox{\o}},V_{\mbox{\o}})$ constructed
above. Clearly $k_1=0$. Let us verify directly that, as predicted
by Proposition~\ref{pisu_v_aute_ve_Dvore},  
the element $\kappa = m + c \in \Lin(\otexp V2,V) \oplus \Lin(V,\otexp
V2)$ satisfies the master equation~(\ref{44}) if and only if
$(V,m,c)$ forms a bialgebra.

From degree reasons, the only possibly nontrivial
components of $h_n(\kappa,\ldots,\kappa)$ are
\[
h_n(\kappa,\ldots,\kappa)^1_3,
h_n(\kappa,\ldots,\kappa)^2_2\ \mbox { and }\
h_n(\kappa,\ldots,\kappa)^3_1.
\]
These values are determined by $\psfb(\jednatri)$, $\psfb(\dvadva)$ and
$\psfb(\trijedna)$. Looking at these values we see that
$h_n(\kappa,\ldots,\kappa) \not= 0$ only for $n=2$ or $n=4$.
The components of $h_2(\kappa,\kappa)$ are
\[
h_2(\kappa,\kappa)^1_3 = 2 \left(\hskip -.3em
\raisebox{-.4em}{\bZbbZdec mm} -
\raisebox{-.4em}{\ZbbZbdec mm} \hskip .2em
\right)\hskip -.2em  ,\
h_2(\kappa,\kappa)^1_3 = 2 \left(\hskip -.3em
\raisebox{-.4em}{\bVbbVdec mm} -
\raisebox{-.4em}{\VbbVbdec mm} \hskip .2em
\right)\hskip -.2em ,\
h_2(\kappa,\kappa)^2_2 = 2 \raisebox{-.4em}{\dvojiteypsilondec cm}
\hskip -.2em .
\]
The components of $h_4(\kappa,\kappa,\kappa,\kappa)$ are
\[
h_4(\kappa,\kappa,\kappa,\kappa)^1_3 =
h_4(\kappa,\kappa,\kappa,\kappa)^3_1 = 0\
\mbox { and }\
h_4(\kappa,\kappa,\kappa,\kappa)^2_2 = 24 \hskip .3em
\frac{\jednadvadec m  \jednadvadec m}
     {\rule{0em}{1.8em}\dvajednadec c \dvajednadec c} \hskip .2em.
\]
The above calculations make the claim obvious.
}
\end{example}

\section*{Appendix~A: Symmetric brace algebras}
\label{apA}

The material of this appendix is taken almost word-by-word
from~\cite{lm:sb}.

\begin{definition}
A symmetric brace algebra is a graded vector space $W$ together with a
collection of degree $0$ multilinear braces $x\langle
x_1,\ldots,x_n\rangle$ that are graded symmetric in $x_1,\ldots,x_n$
and satisfy the identities $ x\langle \hskip .5em \rangle=x$ and
\begin{eqnarray}
\label{axiom}
\lefteqn{\hskip -1em
x\langle x_1,\ldots,x_m\rangle\langle y_1,\ldots,y_n\rangle=}
\\
\nonumber
&&
\sum
\epsilon \cdot x\langle x_1\langle y_{i_1^1},\sqldots,
y_{i_{t_1}^1}\rangle,
x_2\langle y_{i_1^2},\sqldots,y_{i_{t_2}^2}\rangle,\sqldots,
x_m\langle y_{i_1^m},\sqldots,y_{i_{t_m}^m}\rangle,y_{i_1^{m+1}},\sqldots,
y_{i_{t_{m+1}}^{m+1}}\rangle
\end{eqnarray}
where the sum is taken over all unshuffle decompositions
$$
i_1^1<\cdots<i_{t_1}^1,\ldots,i_1^{m+1}<\cdots<i_{t_{m+1}}^{m+1}
$$
of $\{1,\ldots,n\}$ and where $\epsilon$ is the Koszul sign of the
permutation
\begin{eqnarray*}
\lefteqn{
(x_1,\ldots,x_m, y_1,\ldots,y_n) \longmapsto \hskip 3em}
\\
&&  \hskip 3em
(x_1,  y_{i_1^1},\ldots, y_{i_{t_1}^1},
x_2,y_{i_1^2},\ldots,y_{i_{t_2}^2},\ldots,
x_m,
y_{i_1^m},\ldots,y_{i_{t_m}^m},y_{i_1^{m+1}},\ldots,
y_{i_{t_{m+1}}^{m+1}})
\end{eqnarray*}
of elements of $W$.
\end{definition}

For elements $x,y$ of an arbitrary symmetric brace algebra $W$, put $x
\di y := x \angles y$.  One easily proves that then $(W,\di)$ is a
graded pre-Lie algebra in the sense
of~\cite[Section~2]{gerstenhaber:AM63}.

Vice versa, higher brackets $x\langle \Rada x1n
\rangle$ of an arbitrary symmetric brace algebra are,
for $n \geq 2$, determined by their `pre-Lie part' $x \di y = x\langle
y \rangle$. For
instance, axiom~(\ref{axiom})
implies that $x \langle x_1,x_2 \rangle$ can be expressed as
\[
x \langle x_1,x_2 \rangle =  x \sb {x_1} \sb {x_2} -
x\sb {x_1 \sb {x_2}} = (x \di x_1) \di x_2 - x \di (x_1
\di x_2).
\]
The same axiom applied on $x\sb{\Rada x1{n-1}}\sb{x_n}$ can then be
clearly
interpreted as an inductive rule defining $x\sb{ \Rada x1n}$
in terms of $x\langle \Rada x1k \rangle$, with $k < n$.

As proved in~\cite{guin-oudom}, {\em an arbitrary\/} pre-Lie algebra
determines in this way a {\em unique\/}
symmetric brace algebra. Let us emphasize that
this statement is not obvious. First,
axiom~(\ref{axiom}) interpreted as an
inductive rule is `overdetermined.' For example, $x\sb {x_1,x_2,x_3}$
can be expressed both from~(\ref{axiom})
applied to $x \sb {x_1,x_2}\sb {x_3}$ and also
from~(\ref{axiom}) applied to $x \sb {x_1}
\sb {x_2,x_3}$, and it is not obvious whether the results would be the
same. Second, even if the braces are well-defined, it is not clear
whether they satisfy the axioms of brace algebras, including the
graded symmetry.

\section*{Appendix B: 
          Proof of formula~(\protect\ref{Merkulov_je_rusky_hovado})}
\label{apB}

In this appendix we indicate how to prove that the Gerstenhaber-Schack
differential is the same as the differential given by
formula~(\ref{je_to_takova_glegla}) that involves the minimal model
$\sfM$ of the bialgebra \PROP\ $\sfB$. This can in principle also be
achieved by
analyzing the explicit minimal model described
in~\cite{umble-saneblidze:KK}, but we sketch out a procedure that
requires much less information. We will
explain it for the case of bialgebras, but it will be clear how to
generalize it to algebras over an arbitrary \PROP.

Let us start by rewriting~(\ref{je_to_takova_glegla}) in a way that
will make it clear which part of the minimal model it really needs. 
As we already remarked, we know that the minimal model of the
bialgebra \PROP\ $\sfM$ is of the form
$(\freePROP(\{\xi^m_n\}),\psfb)$, 
where $\xi^m_n$ are, for $(m,n) \in I$ defined in~(\ref{SM_je_hovado}),  
generators of biarity $(m,n)$ and degree $m+n -3$.  

Consider the free \PROP\
\[
\sfDM :=
\freePROP(\{\xi^m_n\},\{\eta^m_n\},\varphi),\  (m,n) \in I,
\]  
where $\xi^m_n$ are the generators of the \PROP\ $\sfM$, $\varphi$ is
a generator of biarity $(1,1)$ placed in degree~$0$, and $\eta^m_n :=
\ \susp \xi^m_n$ are, for $(m,n) \in I$, the `suspended' generators of $\sfM$,
i.e.~$\eta^m_n$ is of biarity $(m,n)$ and of degree $m+n -2$.
Introduce finally the degree $+1$ derivation $\dd : \sfDM \to \sfDM$ by
\[
\dd (\xi^m_n) := \eta^m_n,\ \dd (\eta^m_n):= 0,\ \mbox { for }
(m,n) \in I,\ \mbox { and } \dd(\varphi) := 0. 
\]  
We use $\dd$ to equip $\sfDM$ with the differential  $\pa_D$ 
given by
\begin{eqnarray}
\nonumber 
\pa_D(\xi^m_n) & := & \psfb(\xi^m_n),\ \mbox { for $(m,n) \in I$, ($\psfb$
is the differential in $\sfM$)}
\\
\nonumber 
\pa_D(\varphi) &:=& 0, \mbox { and}
\\
\label{hajzl_ruskej}
\pa_D(\eta^m_n)& :=& \varphi^{[m]} \circ \xi^m_n -  \xi^m_n \circ \varphi^{[n]}
 - \dd(\psfb(\xi^m_n)),\ \mbox { for $(m,n) \in I$.}
\end{eqnarray}
where
\begin{equation}
\label{svine_ruska}
\varphi^{[k]} := \sum_{0 \leq j \leq k-1} (\id^{\ot j} \ot \varphi \ot
\id^{\ot (k-j-1)}),\ k \geq 1,
\end{equation}
with $\id \in \sfDM(1,1)$ the unit and $\circ$ denoting the horizontal
composition. 
To prove that $\pa_D^2 = 0$, one needs to observe that
\begin{equation}
\label{hovado_rusky}
\pa_D(\dd (x)) := \varphi^{[m]} \circ x -  x \circ \varphi^{[n]}
 - \dd(\psfb(x))
\end{equation}
for each $x \in \sfM \subset \sfDM$; notice
that~(\ref{hajzl_ruskej}) is~(\ref{hovado_rusky}) with $x = \xi^m_n$. Then
\begin{eqnarray*}
\pa_D^2(\eta^m_n)& =& 
\pa_D\left(\varphi^{[m]} \circ \xi^m_n -  \xi^m_n \circ \varphi^{[n]}
 - \dd(\psfb(\xi^m_n))\right) 
\\
&=&\varphi^{[m]} \circ \psfb(\xi^m_n) -  \psfb(\xi^m_n) \circ \varphi^{[n]}
 - \psfb(\dd(\psfb(\xi^m_n))) = 0,
\end{eqnarray*}
where the vanishing in the second line follows
from~(\ref{hovado_rusky}) applied to $x = \psfb(\xi^m_n)$. The vanishing
$\pa_D^2(\xi^m_n) =\pa_D^2(\varphi) = 0$ is obvious. 

Representing $\psfb(\xi^m_n)$ as in~(\ref{je_to_takova_glegla}), 
the rightmost term $\dd(\psfb(\xi^m_n))$ in~(\ref{hajzl_ruskej}) can be
written as
\begin{equation}
\label{hajzl}
\dd(\psfb(\xi^m_n)) := 
\sum_{s \in S^m_n} \sum_{v \in \Vert(G_s)} \epsilon^s \cdot \dd{G^v_s},
\end{equation}
with $\dd{G^v_s}$ the decorated graph obtained from $G_s$ by
replacing the vertex $v$ of biarity, say, $(p,q)$ (which is by definition
decorated by $\xi^{p}_{q}$) by the vertex of the same biarity decorated
by~$\eta^{p}_{q}$. Loosely speaking, $\dd{G^v_s}$ is obtained by
raising the degree of $v$ by one.

Let $\rho : \sfM \to \sfB$ be the minimal model map and
$\hat\rho_D$ the composition
\[
\hat\rho_D : \sfDM \cong \freePROP(\{\xi^m_n\})
\copr \freePROP(\{\eta^m_n\},\varphi) \stackrel{\rho \copr
  \iden}{\longrightarrow}  
\sfB \copr \freePROP(\{\eta^m_n\},\varphi).
\]
We equip $\sfB \copr \freePROP(\{\eta^m_n\},\varphi)$ 
with the differential $\bpd$ given by the formulas
\[
\bpd|_\sfB = 0,\ \bpd(\varphi) := 0 \mbox { and }
\bpd(\eta^m_n) := \hat\rho (\pa_D(\eta^m_n)),\ \mbox{ for } (m,n) \in I.
\]

Consider the 
$\sfB$-submodule $\sfB\langle \{\eta^m_n\},\varphi \rangle$ of the
coproduct
$\sfB \copr \freePROP(\{\eta^m_n\},\varphi)$ generated by 
$\{\eta^m_n\}_{(m,n) \in I}$ and
$\varphi$. It is spanned by monomials in $\sfB \copr
\freePROP(\{\eta^m_n\},\varphi)$ containing precisely one $\eta^m_n$
or $\varphi$. It is clear that $\sfB\langle
\{\eta^m_n\},\varphi \rangle$ is $\bpd$-stable and that it in fact represents
the free $\sfB$-module generated by $\{\eta^m_n\}_{(m,n) \in I}$ and
$\varphi$. 

On the other hand, for arbitrary $(p,q)$, the \PROP{ic} structure of
$\End_V$ determines an $\End_V$-module action on the suspension
$\uparrow^{p+q-1} \End_V$. Given a bialgebra $B =
(V,\mu,\Delta)$, the corresponding map $\alpha : \sfB \to \End_V$
composed with this action determines
a $\sfB$-module structure on $\uparrow^{p+q-1} \End_V$.

By the freeness of $\sfB\langle \{\eta^m_n\},\varphi \rangle$, 
each $f \in \CGSbi pq = \Lin(\otexp Vp,\otexp Vp)$ specifies the
$\sfB$-module map  
\[
\omega_f :
\sfB\langle \{\eta^m_n\},\varphi \rangle  \to\ \uparrow^{p+q-1} \hskip
-.3em\End_V
\]
that satisfies
\[
\omega_f(\phi) := 0,\
\omega_f(\eta^p_q) := f\ \mbox { and } 
\omega_f(\eta^m_n)\ :=\ 0\ \mbox{ for } (m,n)
\not= (p,q).  
\]
Using~(\ref{hajzl}), one can rewrite~(\ref{psano_v_Koline}) as
\begin{equation}
\label{67}
(\delta_\sfB f)^m_n = \omega_f(\bpd(\eta^m_n)),\ (m,n) \in I.
\end{equation}
This equality has an important consequence which we formulate as

\begin{proposition}
The Gerstenhaber-Schack differential is
determined by the dg-\PROP\ 
\[
\overline{\sfDM} := (\sfB \copr \freePROP(\{\eta^m_n\},\varphi), \bpd).
\]
\end{proposition}

In the rest of this appendix we show that the dg-\PROP\ 
$\overline{\sfDM}$ can be
described without the knowledge of the minimal model $\sfM$ of the
bialgebra \PROP\ $\sfB$.
We say that a degree~$0$ map
$\vartheta : V \to V$ is a {\em derivation\/} of a bialgebra 
$B = (V,\mu,\Delta)$ if
\[
\vartheta \circ \mu = 
\mu (\vartheta \ot \iden_V + \iden_V \ot \vartheta)\ \mbox { and }\
\Delta \circ \vartheta = (\vartheta \ot \iden_V + \iden_V \ot \vartheta).
\]
Let $\sfDB$ be the \PROP\ describing structures consisting of a
bialgebra $B$ and a derivation $\vartheta$ of~$B$. Consider
the homomorphism
\[
\brd : (\sfB \copr \freePROP(\{\eta^m_n\},\varphi),\bpd) \to (\sfDB,0)
\]
such that $\brd|_\sfB$ coincides with the canonical
inclusion $\sfB \hookrightarrow \sfDB$, $\brd(\eta^m_n) = 0$ for
$(m,n) \in I$ and $\brd(\phi) \in \sfDB(1,1)$ is the generator for the
derivation. Let us check that $\brd$ defined in this way is a dg-map.

The equation $\brd(\bpd (\varphi))= 0$ is clear.
The vanishing of  \hbox{$\brd(\bpd (\eta^m_n))$} for $m+n > 3$ follows from
degree reasons. It remains to show that \hbox{$\brd(\bpd (\eta^1_2))=
\brd(\bpd (\eta^2_1)) = 0$}.
By~(\ref{hajzl_ruskej}),
\[
\brd(\bpd (\eta^1_2)) = \brd(\hat \rho_D (\pa_D(\eta^1_2))) =
\brd(\hat \rho_D (\phi \circ \xi^1_2 - \xi^1_2 \circ \phi^{[2]}))=
\vartheta \circ \mu - \mu \circ (\vartheta \ot \id + \id \ot \vartheta),
\] 
where we denoted by the same symbols operations on $V$ and
the corresponding generators in~$\sfDB$; we are sure this will not
lead to a confusion here.
We conclude that $\brd(\bpd (\eta^1_2))$ vanishes because $\vartheta$ is a
$\mu$-derivation. The vanishing of $\brd(\bpd (\eta^2_1))$
can be proved in the same manner. The central
statement of this section is:

\begin{proposition}
\label{zitra_na_ZS_a_neni_snih}
The map $\brd : (\overline{\sfDM},\bpd) \to (\sfDB,0)$ is a homology 
isomorphism.
\end{proposition}

The above proposition says that $(\overline{\sfDM},\bpd)$ is
{\em a $\sfB$-free minimal model\/} of the  \PROP\
$\sfDB$. By the $\sfB$-freeness we mean that $\overline{\sfDM}$ is
obtained by adding free generators to $\sfB$. Minimality means that
the image of the differential $\bpd$ consists of decomposable elements
of the augmented \PROP\ $\overline{\sfDM}$ --
see~\cite[page~344]{markl:ba} for a definition of indecomposables in
augmented \PROP{s}. 

Assuming uniqueness of minimal models, {\em any\/}
$\sfB$-free minimal model of $\sfDB$ will
be isomorphic to  $\overline{\sfDM}$.  
One candidate for  such a model can be constructed by expanding formulas for
the Gerstenhaber-Schack differential into diagrams; denote the
\sfB-free 
dg-\PROP\ obtained in this way 
by $(\overline{\sfGS},\pa_{{\sf GS}})$. 
Methods developed in~\cite{mv} then can be used
to prove that the canonical map $(\overline{\sfGS},\pa_{{\sf GS}}) 
\to (\sfDB,0)$ is a
homology isomorphism. An analog of~(\ref{67}) for
$(\overline{\sfGS},\pa_{{\sf GS}})$ 
then identifies $\delta_\sfB$ with the Gerstenhaber-Schack differential.

\begin{proof}[Proof of Proposition~\ref{zitra_na_ZS_a_neni_snih}]
Let $\rho_D : \sfDM \to \sfDB$ be the composition $\overline{\rho}_D
\circ \hat{\rho}_D$. The proposition is a combination of the
following two statements:

(i) The map $\rho_D : \sfDM \to \sfDB$ is a homology isomorphism,
    i.e.~$\sfDM$ is a minimal model of~$\sfDB$.

(ii) The map $\hat \rho_D : \sfDM \to \overline{\sfDM}$ 
    is a homology isomorphism, too.

Let us prove~(i) first. Consider the grading $\gr$ of $\sfDM$ defined by 
\[
\gr(\phi) := 1\ \mbox { and }\ \gr(\xi^m_n) = \gr(\eta^m_n) :=
0\ \mbox { for }\ (m,n) \in I, 
\]   
and decompose $\pa_D$ as $\pa_D = \pa_1 + \pa_2$, where $\pa_1$ raises
the grading by one and $\pa_2$ preserves it. One can easily check that
$\pa_1^2 = \pa_2^2 = 0$ and that $\pa_1\pa_2 + \pa_2\pa_1 = 0$,
therefore $(\sfDM,\pa_1+\pa_2)$ is a bicomplex. It is not difficult to
prove that 
\begin{equation}
\label{45}
H_*(\sfDM,\pa_1) \cong \freePROP(\{\xi^m_n\},\varphi)/{\sf I},
\end{equation}
with the \PROP{ic} ideal ${\sf I}$ generated by the relations
\[
\varphi^{[m]}\circ \xi^m_n \circ \varphi^{[n]} = 0,\ (m,n) \in I,
\]
which say that $\varphi$ is a derivation with respect to all $\xi^m_n$,
see~(\ref{svine_ruska}) for the notation. It follows from~(\ref{45})
that the $E^1$-term of the spectral sequence associated to this
bicomplex equals
\[
(E^1,d^1) \cong (\freePROP(\{\xi^m_n\},\varphi)/{\sf I},\pa),
\]
where $\pa$ coincides with the minimal model differential $\pa_\sfB$ on the
$\xi$-generators and $\pa(\varphi) = 0$. We conclude that
$H_*(E^1,d^1) \cong \sfDB$, so the associated spectral sequence
degenerates at the $E^2$-level from degree reasons and
$H_*(\sfDM,\pa_D) \cong \sfDB$. One can easily see that the latter
isomorphism is induced by $\rho_D$. This proves~(i).

To prove~(ii), one introduces the grading $\gr'$ of $\sfDM$ by formulas
\[
\gr'(\xi^m_n) := 0,\ \gr'(\phi) := 2\ \mbox { and }\ \gr'(\eta^m_n) :=
m+n\ \mbox { for }\ (m,n) \in I. 
\] 
and a similar grading $\gr''$ of $\overline{\sfDM}$ by
\[
\gr''(b) := 0\ \mbox { for }\ b \in \sfB,\ 
\gr''(\phi) := 2\ \mbox { and }\ \gr''(\eta^m_n) :=
m+n\ \mbox { for }\ (m,n) \in I.
\]
Since $\hat \rho_D$ preserves these gradings, it
preserves also the filtrations
\[
F'_p := \{x \in \sfDM;\ \gr'(x) \leq p\}\ \mbox { and }\
F''_p := \{x \in \overline{\sfDM};\ \gr''(x) \leq p\}
\]
and induces the map $\{E_p(\hat \rho_D) : (E'_p,d'_p) \to (E''_p,d''_p)\}$ 
of the induced spectral sequences. Clearly
\[
(E'_0,d'_0) \cong (\freePROP(\{\xi^m_n\},\{\eta^m_n\},\varphi),\pa),
\]
where $\pa$ coincides with the differential $\pa_\sfB$ of the minimal model
$\sfM$ of $\sfB$ 
on the $\xi$-generators and is trivial on the remaining ones, 
so $H_*(E'_0,d'_0) \cong \sfB \copr
\freePROP(\{\eta^m_n\},\varphi)$. On the other hand, clearly $d''_0
= 0$, therefore $H_*(E''_0,d''_0) \cong \sfB \copr
\freePROP(\{\eta^m_n\},\varphi)$, too. In fact, it is not hard to
see that 
\[
E_1(\hat \rho_D) : (E'_1,d'_1) \to (E''_1,d''_1)
\] 
is an isomorphism of complexes. A
standard spectral sequence argument then finishes our proof of~(ii).
\end{proof}


\end{document}